\RequirePackage{fix-cm}
\documentclass[smallextended,hidelinks,sort&compress]{svjour3}       
\smartqed  
\usepackage{graphicx}
\usepackage{comment}
\usepackage{amsmath,amssymb}
\usepackage{svg, color, subfig}
\usepackage{siunitx}
\usepackage{pgfplots}
\pgfplotsset{compat=1.16}	
\usetikzlibrary{patterns}

\usepackage{booktabs}

\definecolor{azure}{rgb}{0.0, 0.5, 1.0}
\definecolor{awesome}{rgb}{1.0, 0.13, 0.32}
\definecolor{asparagus}{rgb}{0.53, 0.66, 0.42}
\definecolor{cadetgrey}{rgb}{0.57, 0.64, 0.69}
\definecolor{amber}{rgb}{1.0, 0.49, 0.0}

\DeclareSymbolFont{mathx}{U}{mathx}{m}{n}

\usepackage{hyperref}

\usepackage{printlen}

\newcommand{\K}{\mathbb{K}}

\newcommand{\W}{\mathcal{W}}
\newcommand{\LL}{\mathcal{L}}
\newcommand{\R}{\mathbb{R}}

\newcommand{\N}{\mathbb{N}}

\newcommand{\norm}[1]{\left\lVert#1\right\rVert}
\newcommand{\normiii}[1]{\left\vvvert#1\right\vvvert}
\newcommand{\abs}[1]{\left\lvert#1\right\rvert}

\newtheorem{thm}{Theorem} 

\newtheorem{defn}{Definition}
\newtheorem{rem}{Remark}

\newtheorem{propn}{Proposition}

\newcommand{\ee}{\mathbf{e}}
\newcommand{\xx}{\mathbf{x}}
\newcommand{\yy}{\mathbf{y}}

\newcommand{\rr}{\mathbf{r}}
\newcommand{\ff}{\mathbf{f}}

\newcommand{\bb}{\mathbf{b}}

\newcommand{\uu}{\mathbf{u}}
\newcommand{\UU}{\mathbf{U}}
\newcommand{\vv}{\mathbf{v}}
\newcommand{\ww}{\mathbf{w}}

\newcommand{\phiB}{\pmb{\phi}}

\newcommand{\bx}{\mathbf{X}}

\newcommand{\bu}{\mathbf{U}}

\begin{document}

\title{Quasistatic fracture evolution
}


\author{Debdeep Bhattacharya         \and
        Robert Lipton \and
        Patrick Diehl
}


\institute{D. Bhattacharya \at
             Department of Mathematics, Louisiana State University,  Baton Rouge, LA 70803,  USA \\ 
              \email{debdeepbh@lsu.edu}           
              \and
           R. P. Lipton \at
             Department of Mathematics, Louisiana State University,  Baton Rouge, LA 70803,  USA \\ 
              \email{lipton@lsu.edu}           
              \and
           P. Diehl \at
             Center for Computation and Technology, Louisiana State University,  Baton Rouge, LA 70803,  USA \\ 
              \email{pdiehl@cct.lsu.edu}           
}

\date{Received: date / Accepted: date}

\maketitle

\begin{abstract}
Nonlocal quasistatic fracture evolution for interacting cracks is developed and supporting numerical examples are presented. The approach is implicit and is based on local stationarity and fixed point methods. 
It is proved that the fracture evolution decreases stored elastic energy with each load step as the cracks advance; provided the load increments are chosen sufficiently small. This is also seen in the numerical examples. 
 The numerical examples include evolution of a straight crack, a crack propagating inside an L-shaped domain, and two offset inward propagating cracks.
\end{abstract}

\section{Introduction}
The hallmark of Peridynamic (PD) simulations is that dynamic fracture patterns can emerge from the nonlocal model and are not prescribed \cite{silling2000reformulation,silling2007peridynamic}.
In the absence of inertia one considers quasistatic or rate independent evolution. We address quasistatic fracture evolution using a nonlocal implicit framework. As before, it is possible for fracture patterns to emerge from the model. It is shown explicitly how this can be seen through bond breaking for displacement controlled evolution. In this treatment, a cohesive law is used to model nonlocal forces. It is shown here that the fracture evolution decreases stored elastic energy with each load step in the undamaged material as the cracks advance; provided the load increments are chosen sufficiently small. We provide a rigorous existence theory of quasistatic nonlocal fracture evolution using fixed point methods. We provide numerical examples and to reduce the computational cost of the implicit scheme we apply the analytic method presented in~\cite{diehl2022quasistatic,bhattacharya2021peridynamics}. Here, the work needed to assemble the tangent stiffness matrix is reduced by using explicit analytic formulas for the Hessian.

There is now a  large literature on dynamic simulations using PD, e.g., see the reviews \cite{bobaru2016handbook,javili2019peridynamics,isiet2021review,diehl2019review,diehl2022comparative}. Building on this, dynamic relaxation methods have been applied to quasistatic PD~\cite{kilic2010adaptive,Yaghoobi}. 
However, the literature addressing implicit quasistatic PD simulations is far smaller than for dynamic simulations, and some representative publications  include~\cite{huang2015improved,mikata2012analytical,zaccariotto2015examples,wang2019studies,breitenfeld2014quasi,kilic2010adaptive,rabczuk2017peridynamics,freimanis2017mesh}.  One reason for the paucity of results is that implicit methods incur large computational costs. 
On the other hand, dynamic relaxation for a quasistatic problem can take numerous iterations to converge. Because of this reason, the advantage of implicit methods for quasistatic nonlocal  problems is that convergence is achieved using relatively few iterations \cite{zaccariotto2019}.  Recently, several techniques to reduce the computational costs for quasistatic simulations are available. These include the adaptive use of linear elastic and peridynamic meshes \cite{zaccariotto2019}, the fast convolution method \cite{jafarzadeh2021general}, fast Galerkin methods~\cite{wang2012fast}, a multi-threaded approach for generating sparse stiffness matrices ~\cite{prakash2020multi}, a combined implicit-explicit method~\cite{hu2016bond2,hu2018thermomechanical}, and the fire  algorithm~\cite{shiihara2019FastQN}. 




In this paper, we use the analytic formulas for bond stiffness~\cite{diehl2022quasistatic,bhattacharya2021peridynamics} to compute the discrete Hessian used in the Newton-Raphson scheme to find the elastic displacement field for each load step. Unlike local fracture theories involving an explicit crack and elastic equilibrium equations off the crack, the nonlocal equations of elastic equilibrium are well-defined everywhere.  Away from the crack the nonlocal solutions are close to local solutions of the elastic equilibrium equation. This is provably true for the nonlocal  modeling of cracks in the limit of vanishing nonlocality, see \cite{lipton2014dynamic,lipton2016cohesive,liptonjha2021}.

We apply the nonlocal model to simulate the fracture evolution of a straight mode-I crack inside a square plate,  capture the emergent crack growth at the re-entrant corner of an L-shaped panel as seen in experiment, and model the interaction of two inward propagating cracks  as they approach each other.
Nonlocal simulations show that the energy inside the undamaged material  decreases with the load step while cracks propagate.
In earlier work ~\cite{https://doi.org/10.1002/nme.7005}  this constitutive model 
is used in a simple comparison of quasistatic damage evolution between force and  displacement  loading and is in line with the theory of crack resistance~\cite{anderson2017fracture}.

The paper is structured as follows: Section~\ref{sec:nonlocal:model} introduces the nonlocal model and displacement controlled bond breaking and fracture. In Section~\ref{energy inequality} Continuity and asymptotic energy reduction for displacement controlled fracture evolution is established. The rigorous existence theory of nonlocal fracture evolution is presented in Section~\ref{exist}.
The algorithm and discritization for the implicit method is given In Section~\ref{Newt}. In Section~\ref{sec:numerics} numerical results and comparison with experimental data are presented.



\section{Nonlocal model and displacement controlled bond breaking and fracture}
\label{sec:nonlocal:model}
\begin{figure}
    \centering
\includegraphics[width=0.8\linewidth]{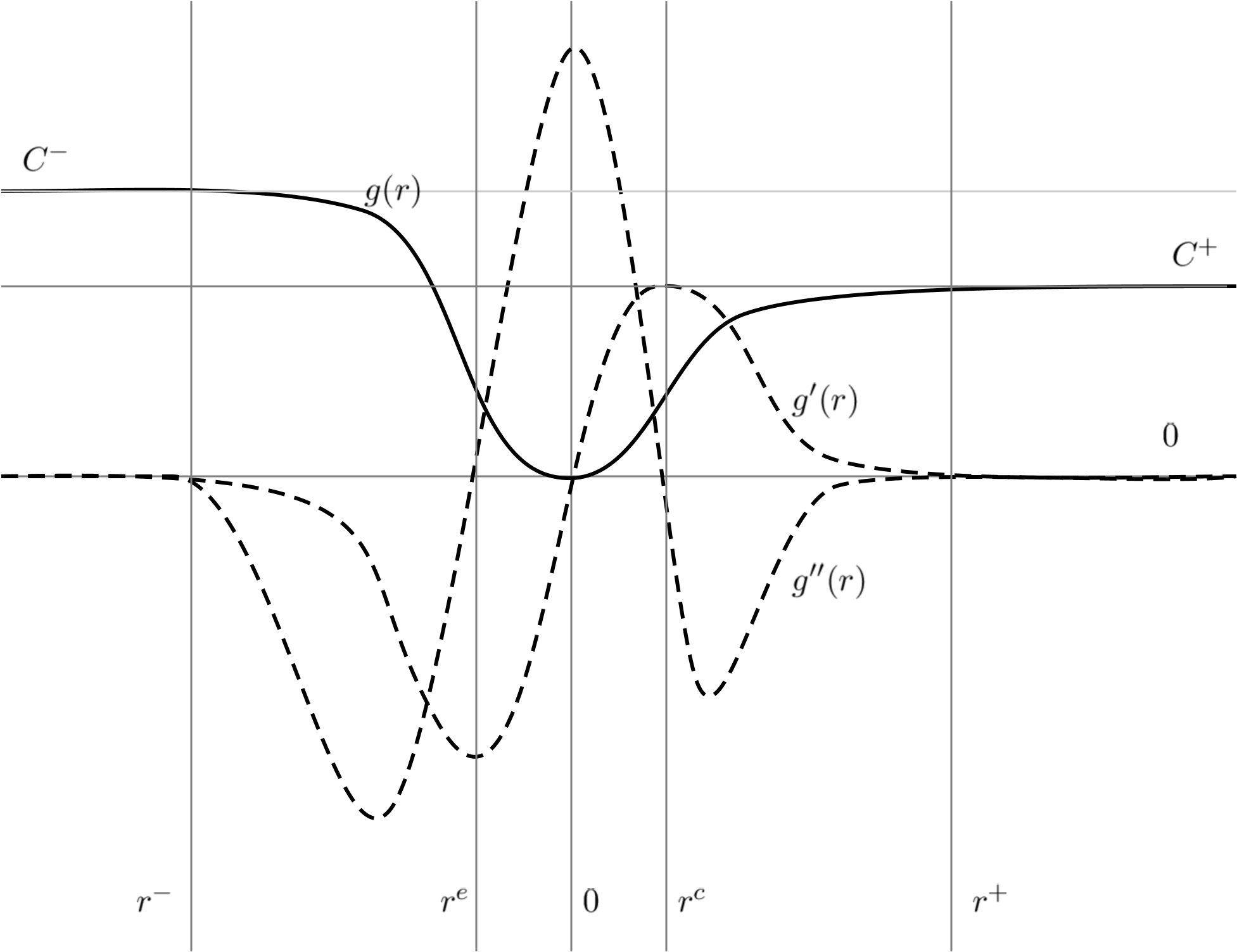}
		  \caption{The potential function $g(r)$ and derivatives $g'(r)$ and $g''(r)$  for tensile force. Here $C^+$ and $C^-$ are the asymptotic values of $g$. The derivative of the force potential goes smoothly to zero at $r^+$ and $r^-$.}
		  \label{ConvexConcavea}
\end{figure}

A displacement controlled fracture evolution is addressed.  We provide a nonlocal mesoscopic model for quasistatic fracture. The Dirichlet data for the prescribed displacement is specified on an interaction domain $\Omega_d$. We introduce $\Omega$, with $\Omega_d\subset \Omega$ and the cracking body $D=\Omega\setminus\Omega_d$. The interaction domain $\Omega_d$ is of thickness equal to the length scale of non-local interaction $\epsilon$. Here $\Omega$ is a bounded domain in two or three dimensions. Nonlocal interactions between a point $\xx$ and its neighbors $\yy$ are confined to the sphere (disk) $H_\epsilon(\xx)=\{\yy:\,|\yy-\xx|<\epsilon\}$. The radius $\epsilon$ is called the horizon and is chosen an order of magnitude smaller than the length scale of the domain $\Omega$.
We introduce the 
nonlocal strain $S(\yy, \xx, \uu)$ between the point $\xx$ and any point $\yy \in H_\epsilon(\xx)$  given by
\begin{align*}
    S(\yy, \xx, \uu) = \frac{\uu(\yy) - \uu(\xx)}{\abs{\yy - \xx}} \cdot \ee_{\yy - \xx},
\end{align*}	
where $\ee_{\yy- \xx}$ is the unit vector given by
\begin{align*}
\ee_{\yy - \xx}  = \frac{\yy - \xx}{\abs{\yy - \xx}}.
\end{align*}	
Force is related to strain using the constitutive relation given by the cohesive force law as in \cite{lipton2014dynamic,lipton2016cohesive}. Under this law the force is linear for small strains and for larger strains the force begins to soften and then approaches zero
after reaching a critical strain.  The force function is $g'$ is shown in figure \ref{ConvexConcavea}.
The nonlocal force density $\ff$ is given in terms of the nonlocal potential $\W(S)$ by
\begin{align}\label{contsit1}
    \ff(\yy, \xx, \uu) = 2 \partial_S \W(S(\yy, \xx, \uu)) \ee_{\yy - \xx},
\end{align}	
where
\begin{align}\label{contsit2}
    \W(S(\yy, \xx, \uu)) = \frac{J^\epsilon(\abs{\yy - \xx})}{\epsilon^{n+1}\omega_n\abs{\yy - \xx}}  g(\sqrt{ \abs{\yy - \xx}} S(\yy, \xx, \uu)).
\end{align}	
Here, $J^\epsilon(r) = J(\frac{r}{\epsilon})$, where $J$ is a non-negative bounded function supported on $[0,1]$. $J$ is called the \textit{influence function} as it determines the influence of the bond force of peridynamic neighbors $\yy$ on the center $\xx$ of $H_\epsilon(\xx)$. The volume of unit ball in $\R^n$ is denoted by $\omega_n$. 
As figure \ref{ConvexConcavea} illustrates, we assume that $g(r)$ and the derivatives $g'(r)$, $g''(r)$, and $g'''(r)$ are bounded for $-\infty < r<\infty$. It is required is that $g(0)=0$ and $g(r)>0$ otherwise, $g(r)$ together with its first three derivatives must be bounded, and that $g$ be  convex in the interval $r^e<0<r^c$ and concave outside this interval with finite limits $\lim_{r\rightarrow-\infty}{g(r)}=C^-$ and $\lim_{r\rightarrow\infty}{g(r)}=C^+$. Additionally $\max\{|g''(r|)\}=g''(0)$.

The influence function and $g$ are calibrated for a given material with known Lam\'e modulus $\mu$ and critical energy release rate $G_c$ using the relations

\begin{align*}
    \mu  = g''(0)/ 10 \int_0^1 r^3 J(r) dr, \ n = 3,
    \\
    \mu  = g''(0)/ 8 \int_0^1 r^2 J(r) dr, \ n = 2,
\end{align*}
and
\begin{align*}
    G_c = 2 \frac{\omega_n}{\omega_{n-1}} g_\infty \int_0^1 r^n J(r) dr, \ n=2,\,3.
\end{align*}
Here we are using bond based nonlocal interactions so for plane stress $\nu=1/3$ and in three dimensions $\nu=1/4$ and $\lambda=\mu$ in both cases. We can also apply more general force interactions as in state based nonlocal interactions for softening models, see \cite{Lipton2018free}.

\subsection{Bond breaking evolution}
\label{bomdbreaking}
It is clear from figure \ref{ConvexConcavea} that $g''(r)<0$ for
\begin{align}\label{broken}
    {r^e}>{r=\sqrt{|\yy - \xx|}}S(\yy, \xx, \uu) \hbox{ or }
  \sqrt{|\yy-\xx|} S(\yy, \xx, \uu)&= r > {r^c},
\end{align}
and we say that the bond is broken between $\yy$ and $\xx$ when 
\begin{align} \label{breakingbonds}
g'(\sqrt{|\yy-\xx|}S(\yy,\xx,\uu))=0 \hbox{  and the strain satisfies \eqref{broken}.}
\end{align}
Initially the characteristic function for all intact bonds between points  $\yy\in\Omega$ inside  $H_\epsilon(\xx)$ and $\xx\in D$  is denoted by $\chi_0(\yy,\xx)$ and initially  the whole domain is intact so  $\chi_0(\yy,\xx)=1$.

The nonlocal force density $\LL_0$ is defined for all points  $\xx$ in $D$ and $\yy$ in $\Omega$ is given by
\begin{align}
      \label{eq: force}
&  \LL_0 [\uu](\xx) = \nonumber\\
&-\int\limits_{H_\epsilon(\xx) \cap \Omega}^{} {2}\chi_0(\yy,\xx)\frac{J^\epsilon(\abs{\yy - \xx} )}{\epsilon^{n+1} \omega_n \sqrt{\abs{\yy - \xx}}} g'\left(\sqrt{ \abs{\yy - \xx}} S(\yy, \xx, \uu)\right) \ee_{\yy - \xx} d\yy.
\end{align}

The solution to the nonlocal boundary value problem is given by a function $\uu^1=\UU^1$ on $\Omega_d$ for which
\begin{equation}\label{eq:bvp}
\LL _0[\uu^1](\xx) =0,\,\,\hbox{  for  $\xx$ in $D$}.
\end{equation}

Given the solution, $\uu^1$ consider all pairs $(\yy,\xx)$  with  $\yy\in\Omega$ inside  $H_\epsilon(\xx)$ and $\xx\in D$  for which the bond between them is broken. 
This set of pairs is called the set  $\Delta S_1$. Denote the new operator $\LL_1(\uu)$ obtained from \eqref{eq: force} by deleting all bond pairs in $\Delta S_1$, let $\chi_1(\yy,\xx)$  denote the indicator function of unbroken bonds and the solution $\uu^1$ of \eqref{eq:bvp} is also  a solution of
\begin{equation}\label{eq:bvp1}
\LL_1 [\uu^1](\xx) =0,\,\,\hbox{  for  $\xx$ in $D$},
\end{equation}
where
\begin{align}
      \label{eq: force1bis}
&  \LL_1 [\uu](\xx) = \nonumber\\
&-\int\limits_{H_\epsilon(\xx) \cap \Omega}^{} {2}\chi_1(\yy,\xx)\frac{J^\epsilon(\abs{\yy - \xx} )}{\epsilon^{n+1} \omega_n\sqrt{\abs{\yy - \xx}}} g'\left(\sqrt{ \abs{\yy - \xx}} S(\yy, \xx, \uu)\right) \ee_{\yy - \xx} d\yy.
\end{align}

Next we increment the boundary load to $\UU^2$ to get the solution $\uu^2$ of 
\begin{equation}\label{eq:bvp21}
\LL_1 [\uu^2](\xx) =0,\,\,\hbox{  for  $\xx$ in $D$}.
\end{equation}
where $\uu=\UU^2$ on $\Omega_d$.
Again given $\uu^2$ consider all pairs $(\yy,\xx)$ of points in $\Omega\times D\setminus\ \Delta S_1$ for which $|\yy-\xx|<\epsilon$ and the bond is broken. 
Call this set of pairs $\Delta S_2$.  Now set $S_1=\Delta S_1$, and set $S_2=S_1\cup \Delta S_2$. Denote the new operator $\LL_2(\uu)$ obtained from \eqref{eq: force} by deleting bond pairs $S_2$ and note that the solution $\uu^2$ of \eqref{eq:bvp21} is also a solution of
\begin{equation}\label{eq:bvp22}
\LL_2 [\uu^2](\xx) =0,\,\,\hbox{  for  $\xx$ in $D$}.
\end{equation}
We can iterate this process with $M$ displacement loads $\{\UU^N\}_{N=1}^M$ and set $S_N=S_{N-1}\cup \Delta S_{N}$ to get a sequence of operators 
\begin{align}
      \label{eq: forcedosN}
  &  \LL_N[\uu](\xx) =\nonumber\\
  & -\int\limits_{H_\epsilon(\xx) \cap \Omega}^{} {2}\chi_N(\yy,\xx)\frac{J^\epsilon(\abs{\yy - \xx} )}{\epsilon^{n+1} \omega_n \sqrt{\abs{\yy - \xx}}} g'\left(\sqrt{ \abs{\yy - \xx}} S(\yy, \xx, \uu)\right) \ee_{\yy - \xx} d\yy,
\end{align}	
solutions $\{\uu^N\}_{N=1}^M$, and debonding sets $\{S_N\}_{N=1}^M$,
with $S_1\subset S_2\cdots \subset S_M$. This constitutes the displacement controlled bond breaking evolution for both two and three-dimensional problems.

\begin{figure} 
\centering
\begin{tikzpicture}[xscale=0.50,yscale=0.50]

\draw [fill=gray, gray] (-8,-0.5) rectangle (-6.0,0.5);





\node [above] at (-2.5,1.5) {Crack after $1^{st}$ load step $C_1=\Delta C_1$};

\node [above] at (-6.5,0.5) {$C_1$};

\draw[fill=gray, gray] (-6.0,-0.5) arc (-90:90:0.5)  -- cycle;












\draw [-,thick] (-8,0) -- (-6.0,0);




\draw [thick] (-8,-5) rectangle (8,5);












\node [above] at (-7.6,0) {$\epsilon$};

\node [below] at (-7.6,0) {$\epsilon$};

\draw [thick] (-7.3,.5) -- (-7.1,.5);

\draw [<->,thick] (-7.215,.5) -- (-7.215,0);

\draw [<->,thick] (-7.215,-.5) -- (-7.215,0);

\draw [thick] (-7.3,-.5) -- (-7.1,-.5);

\draw [->,thick] (-6.0,-1.0) -- (-6.0,-0.1);

\node  at (-6.0,-1.5) { $\ell_1$};



\end{tikzpicture} 

\begin{tikzpicture}[xscale=0.50,yscale=0.50]
\draw [fill=gray, gray] (-8,-0.5) rectangle (2.0,0.5);
\node [above] at (-2.5,0.5) {Crack after $N^{th}$ {\rm load step}:       $ C_N$};
\draw[fill=gray, gray] (2.0,-0.5) arc (-90:90:0.5)  -- cycle;


\draw [-,thick] (-8,0) -- (2.0,0);
\draw [thick] (-8,-5) rectangle (8,5);







\node [above] at (-1.6,0) {$\epsilon$};
\node [below] at (-1.6,0) {$\epsilon$};

\draw [thick] (-1.3,.5) -- (-1.1,.5);

\draw [<->,thick] (-1.215,.5) -- (-1.215,0);

\draw [<->,thick] (-1.215,-.5) -- (-1.215,0);

\draw [thick] (-1.3,-.5) -- (-1.1,-.5);

\draw [->,thick] (2.0,-1.0) -- (2.0,-0.1);

\node  at (2.0,-1.5) { $\ell_N$};



\end{tikzpicture} 
\caption{{ \bf The internal boundary  $C_1=\Delta C_1$ after $1^{st}$ load step given by the line segment of length $\ell_1$ after $N$ steps it grows to length $\ell_N$. The set of broken bonds for the non-local model are given by the grey regions. The set $D_{0}=D$ is the interior of the rectangle, the set $D_{N}$ is obtained by removing $C_N$ at step $N$, see Section \ref{Fracture}.}}
 \label{Crack}
\end{figure}
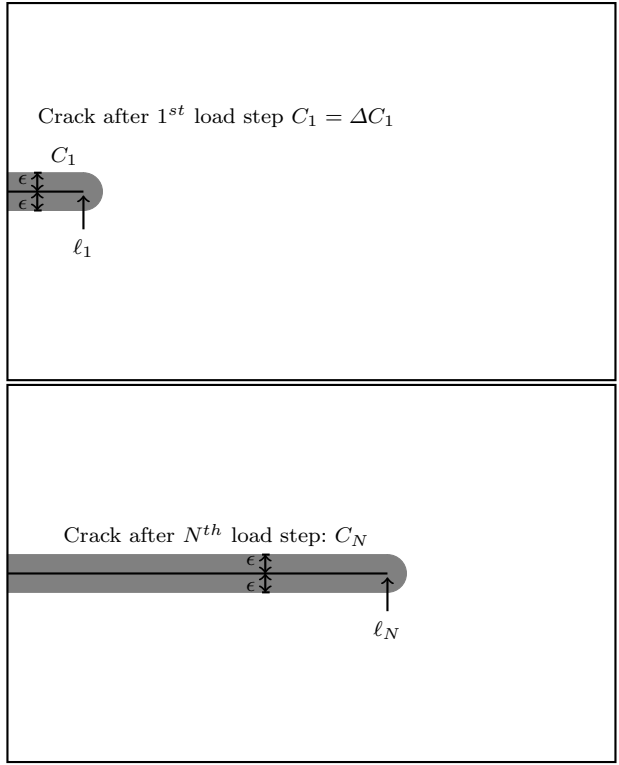

\subsection{Fracture evolution}
\label{Fracture}

It is crucial to observe that crack evolution is not prescribed and only bonds between points are allowed to break. However, a  bond breaking evolution can emerge that is a nonlocal fracture evolution. 
We now illustrate how a bond breaking evolution be can be understood as a nonlocal  fracture evolution.
In what follows we focus on $S_1$ and illustrate the ideas for the 2-dimensional problem since the three-dimensional problem is similar. Suppose the geometry of $S_1$ is characterized by a straight line segment $\Delta C_1$ of length $\ell_1$ across which all bonds of length less than $\epsilon$ are broken. 
For this case $S_1$ corresponds to the union of all neighborhoods that intersect $\Delta C_1$. This is illustrated by the gray region in Figure \ref{Crack}. We delete $\Delta C_1$ from $D$ to form the cracked region $D_1$ with new internal boundary $C_1=\Delta C_1$.  This kind of bond breaking geometry is identical to introducing a new internal boundary $C_1$ associated with a nonlocal traction free boundary condition. {\em In this way a crack appears as an internal boundary $C_1$ with  zero nonlocal traction forces acting on either side.}
So for this case, we can rewrite the operator $\LL_1[\uu](\xx)$ for $\xx$ on $D$ as
\begin{align}
      \label{eq: force1}
    \LL_1 [\uu](\xx) = -\int\limits_{\hat{H}_\epsilon(\xx) \cap \Omega_1} {2}\frac{J^\epsilon(\abs{\yy - \xx} )}{\epsilon^{n+1} \omega_n \sqrt{\abs{\yy - \xx}}} g'\left(\sqrt{ \abs{\yy - \xx}} S(\yy, \xx, \uu)\right) \ee_{\yy - \xx} d\yy,
\end{align}	
where $\Omega_1=\Omega_d\cup D_1$ and $\uu=\UU^1$ on $\Omega_d$. Here, $\hat{H}_\epsilon(\xx)$ referrers to all $\yy$ in $|\yy-\xx|<\epsilon$ located on the same side of the boundary of $D_1$
as $\xx$.  The nonolocal traction free boundary conditions are natural boundary conditions. They are given in \eqref{eq: force1} through the choice of the integration domain $\hat{H}_\epsilon(\xx)\cap\Omega_1$.  More generally, one has a similar formula for $\LL_1$ for the case when $S_1$ given by the collection of all bonds that intersect a smooth curve segment $\Delta C_1$. As before, we delete the smooth curve from $D$ to get $D_1$ and the formula for $\LL_1$ given by \eqref{eq: force1}. Similarly, we can suppose the same for $\Delta S_2,\ldots,\Delta S_M$ to get a growing crack $C_N=C_{N-1}\cup \Delta C_N$,  with decreasing sets  $D_N\subset D_{N-1}\cdots\subset D_1$, $D_N=D_{N-1}\setminus \Delta C_N$ and
operators
\begin{align}
      \label{eq: forceN}
    \LL_N[\uu](\xx) = -\int\limits_{\hat{H}_\epsilon(\xx) \cap \Omega_N} {2}\frac{J^\epsilon(\abs{\yy - \xx} )}{\epsilon^{n+1} \omega_n \sqrt{\abs{\yy - \xx}}} g'\left(\sqrt{ \abs{\yy - \xx}} S(\yy, \xx, \uu)\right) \ee_{\yy - \xx} d\yy,
\end{align}	
with $\Omega_N=\Omega_d\cup D_N$. Here $D_N$ has internal boundary $C_N$ and $\hat{H}_\epsilon(\xx)$ referrers to all pairs $\yy$ in $|\yy-\xx|<\epsilon$ located on the same side of the boundary of $D_N$ as $\xx$.  
This constitutes the displacement controlled fracture evolution. 

The solution $\uu^N$ to $\LL_N[\uu^N]=0$ for each loading $\UU^N$  in the fracture evolution is found numerically using the Newton-Raphson method. This is addressed in section \ref{Newt}. The next section  shows that the elastic energy stored in the intact material decreases with each load step when the load step is sufficiently small.


\section{Continuity and asymptotic energy reduction for displacement controlled fracture evolution}
\label{energy inequality}


Here, it is shown that the fracture evolution decreases stored elastic energy with each load step as the cracks advance; provided the load increments are chosen sufficiently small. This is illustrated when we work on the space of essentially bounded displacements $\uu$ on $\Omega$  denoted by $L^\infty(\Omega,\mathbb{R}^n)$, $n=2,\,3$.
It is assumed that any prescribed displacement $\UU$ on $\Omega_d$ belongs to this space and takes the value $0$ in $D$. This set of boundary displacements is a subspace of  $L^\infty(\Omega,\mathbb{R}^n)$ and we denote it by $\mathcal{B}$. Any displacement $\uu$ where $\uu=\UU$ in $\Omega_d$ can be written as $\uu=\vv+\UU$ where $\vv$ belongs to $\mathcal{V}=\{\vv \in L^{\infty}(\Omega,\mathbb{R}^n),\,\, \vv=0\hbox{ on $\Omega_d$}\}$ and $\UU$ belongs to $\mathcal{B}$.   Adopting standard convention one has $\uu$  in $\mathcal{V}+\mathcal{B}$. We extend $\LL_N[\uu]$ by $0$ for all $\xx$ in $\Omega\setminus D$.

The operator $\LL_N$ satisfies several properties \cite{BhattacharyaLipton}:\\
\noindent The operator $\LL_N$ is uniformly Lipschitz continuous on $L^\infty(\Omega,\mathbb{R}^n)$, i.e.,
\begin{align}\label{lip}
\norm{\LL_N[\uu+\Delta \uu]-\LL_N[\uu]}_\infty\le C\norm{\Delta\uu}_\infty  \hbox{ for $C$ independent of $\uu+\Delta\uu$},
\end{align}
here $||\cdot||_\infty$ is the $L^\infty({\Omega},\mathbb{R}^n)$ norm.  It now follows that $\LL_N[\uu]$  belongs to $\mathcal{V}$ and for fixed $\UU$ in $\mathcal{B}$ the operator $\LL_N[\vv]$ is a bounded operator from $\mathcal{V}$ into itself.

The operator $\LL_N[\uu]$ is  Fr\'echet differentiable and 
    is the bounded linear functional acting  on $\Delta \uu\in L^\infty(\Omega,\mathbb{R}^n)$ given by
\begin{align}
    \label{eq:derivative}
 &   {\LL}'_N[\uu] \Delta \uu = \nonumber\\
 &   - \int\limits_{\hat{H}_\epsilon(\xx) \cap {\Omega}}^{} \frac{J^\epsilon(\abs{\yy - \xx}) }{\epsilon^{n+1} \omega_n}  g''\left( \sqrt{\abs{\yy - \xx}} S\left(\yy, \xx, \uu\right)\right) S(\yy, \xx, \Delta \uu) \ee_{\yy - \xx} d\yy.
\end{align}	
The operator $\LL_N[\uu]$ is continuously Fr\'echet differentiable, i.e.,
\begin{align}\label{freshet}
    \lim_{\Vert\Delta\uu\Vert_\infty\rightarrow 0}\frac{\Vert\LL_N[\uu + \Delta \uu] (\xx) - \LL_N[\uu](\xx) -\LL'_N[\uu]\Delta\uu \Vert_\infty}{\Vert\Delta\uu\Vert_\infty}=0,
\end{align}
and the derivative is Lipshitz continuous in $\uu$, i.e., for $\boldsymbol{\delta}\hbox{, } \Delta\uu \in L^\infty(\Omega,\mathbb{R}^n)$  there is a constant $C$ independent of $\boldsymbol{\delta}$ such that
\begin{align}\label{cont}
    \frac{\Vert \LL'_N[\uu+\boldsymbol{\delta}]\Delta\uu-\LL'_N[\uu]\Delta\uu\Vert_\infty}{\Vert\Delta\uu\Vert_\infty}\leq C\Vert\boldsymbol{\delta}\Vert_\infty. 
\end{align}

From \eqref{cont} we have a $K$ independent of $\uu-\ww$ in $\mathcal{B}+\mathcal{V}$ such that
\begin{align}\label{Lipschits}
    |||\LL'_N[\uu]-\LL'_N[\ww]||| &<  K||\uu-\ww||_\infty
\end{align}
Here $|||\cdot|||$ is the operator norm for the linear functionals defined on $L^\infty(\Omega,\mathbb{R}^n)$.

\begin{rem}[Bond breaking evolution]
\label{bondevo}
In summary  the operator $\LL'_N$ for $\uu$ in $\mathcal{V}+\mathcal{B}$ is given for  both bond breaking and nonlocal fracture by
\begin{align}
      \label{eq: hessianN}
   & \LL'_N[\uu]\Delta\uu(\xx) = \nonumber\\
  &  -\int\limits_{\hat{H}_\epsilon(\xx) \cap \Omega} {\chi_N(\yy,\xx)}\frac{J^\epsilon(\abs{\yy - \xx} )}{\epsilon^{n+1} \omega_n} g''\left(\sqrt{ \abs{\yy - \xx}} S(\yy, \xx, \uu)\right)S(\yy, \xx, \Delta \uu)  \ee_{\yy - \xx} d\yy,
\end{align}	
and properties \eqref{freshet} through \eqref{Lipschits} hold.
\end{rem}

The crack domain at the $N^{th}$ load step is $C_N$ and the set of broken bonds is $S_N$. The set of broken bonds $S_N$ is given by the grey zone in Figure \ref{Crack}. We define the set of intact material by $\tilde{D}_N=D\setminus S_N$. The elastic energy density stored  at a point $\xx$ in $\tilde{D}_N$ is given by
\begin{align*}
  &  W_N(\xx,\uu^N)=\int_{\hat{H}_\epsilon(\xx)\cap\Omega}\,|\yy-\xx|\,\mathcal{W}(\yy,\xx,S(\yy,\xx,\uu^N))\,d\yy.
\end{align*}
The indicator function of the set $\tilde{D}_N$ taking the value $1$ inside $\tilde{D}_N$ and $0$ elsewhere is denoted by $\chi_{\tilde{D}_N}(\xx)$.
The elastic energy of displacement inside the intact material  $\tilde{D}_N$ for prescribed boundary displacement $\UU^N$ on $\Omega_d$ is given by 
\begin{align}\label{intactenergy}
    E_N=\int_{D_0}\chi_{\tilde{D}_N}(\xx)\,W_N(\xx,\uu^N)\,d\xx.
\end{align}

{\bf Energy inequality inside undamaged material}\\
\\We suppose that $\LL'_N[\uu^N]^{-1}$ exists and is bounded for each solution $\uu^N$ in the fracture evolution. When there is crack propagation, i.e.,  $\chi_{\tilde{D}_N}>\chi_{\tilde{D}_{N+1}}$,
the elastic energy satisfies,
\begin{align}\label{equ: energy decrease}
    E_N> E_{N+1}+\omega_{N+1},
\end{align}
If no crack propagation $\chi_{\tilde{D}_N}=\chi_{\tilde{D}_{N+1}}$, then
\begin{align}\label{equ: energy ident}
    E_N= E_{N+1}+\omega_{N+1},
\end{align}
In both cases
\begin{align}\label{equ: limiting1}
    \omega_{N+1}\rightarrow 0, \hbox{ as } ||\UU^{N+1}-\UU^N||_\infty\rightarrow 0.
\end{align}
{\em This shows that the fracture evolution decreases stored elastic energy with each load step as the cracks advance; provided the load increments are chosen sufficiently small.}

One can use these arguments to recover an asymptotic statement about energy reduction.
At the completion of load step $N$  consider a family of load increments $\Delta\UU_\ell$, $\ell=1,2,\ldots$ with $||\Delta\UU_\ell||_\infty\rightarrow 0$ that can be be applied as the $N+1$ load step. 
Let $E_\ell$ to be the elastic energy of the undamaged material corresponding the applied load at the $N+1$ load step is  $\UU^N+\Delta\UU_\ell$, then \eqref{equ: energy decrease}, \eqref{equ: energy ident}, and \eqref{equ: limiting1} imply\\

{\bf Asymptotic energy reduction}
\begin{align}\label{eq: upperfromabove}
E_N\geq\limsup_{||\Delta\UU_\ell||_\infty\rightarrow 0}{E_\ell}.
\end{align}

We establish \eqref{equ: energy decrease} and \eqref{equ: limiting1} noting that \eqref{equ: energy ident} follows. We write the difference as
\begin{align*}
    E_N- E_{N+1}&=\int_{D_0}[\chi_{\tilde{D}_N}(\xx)-\chi_{\tilde{D}_{N+1}}(\xx)]\,W_N(\xx,\uu^N)\,d\xx\nonumber\\
    &+\int_{D_0}\chi_{\tilde{D}_{N+1}}(\xx)\,[W_N(\xx,\uu^N)-W_{N+1}(\xx,\uu^{N+1})]\,d\xx,
\end{align*}
the first term is positive since $\chi_{\tilde{D}_N}(\xx)>\chi_{\tilde{D}_{N+1}}(\xx)$ and $W_N(\xx,\uu^N) \geq 0$. The second term can be written as
\begin{align*}
    \int_{D_0}\chi_{\tilde{D}_{N+1}}(\xx)\,[W_N(\xx,\uu^N)-W_{N+1}(\xx,\uu^{N+1})]\,d\xx= I_1+I_2,
\end{align*}
where $I_1$ and $I_2$ are given by
\begin{align*}
    &I_1=\\
    &\int_{D_0}\chi_{\tilde{D}_{N+1}}(\xx)\,\left(\int_{\hat{H}_\epsilon(\xx)\cap(\Delta D\cup\Omega_d)}\frac{J^\epsilon(\abs{\yy - \xx})}{\epsilon^{n+1}\omega_n\abs{\yy - \xx}}  g(\sqrt{ \abs{\yy - \xx}} S(\yy, \xx, \uu^N))\,d\yy\right)\,d\xx,
\end{align*}
where $\Delta D=D_N\setminus D_{N+1}$ and
\begin{align*}
    I_2&=\int_{D_0}\chi_{\tilde{D}_{N+1}}(\xx)\,(\int_{\hat{H}_\epsilon(\xx)\cap(D_{N+1}\cup\Omega_d)}\frac{J^\epsilon(\abs{\yy - \xx})}{\epsilon^{n+1}\omega_n \abs{\yy - \xx}}  [g(\sqrt{ \abs{\yy - \xx}} S(\yy, \xx, \uu^N))\\
    &-g(\sqrt{ \abs{\yy - \xx}} S(\yy, \xx, \uu^N+\Delta\uu))]\,d\yy)\,d\xx,
\end{align*}
where $\Delta\uu=\uu^{N+1}-\uu^N$. Now $I_1$ is positive since the integrands are positive so
\begin{align*}
    E_N > E_{N+1}+I_2.
\end{align*}
We now set $\omega_{N+1}=I_2$ and get an upper estimate. Using the fundamental theorem of calculus, we have
\begin{align*}
    \omega_{N+1}&=\int_{D_0}\chi_{\tilde{D}_{N+1}}(\xx) \int_0^1T(\xx,t)\,dt\,d\xx,
\end{align*}
where
\begin{align*}
  &  T(\xx,t)=\\
  &\int_{\hat{H}_\epsilon(\xx)\cap(D_{N+1}\cup\Omega_d)}\frac{J^\epsilon(\abs{\yy - \xx})}{\epsilon^{n+1}\omega_n\abs{\yy - \xx}}  g'(\sqrt{ \abs{\yy - \xx}} \rr(t))\sqrt{|\yy-\xx|}S(\yy,\xx,\Delta\uu)\,d\yy,
\end{align*}
with $\rr(t)=\sqrt{|\yy-\xx|}S(\yy,\xx,\uu^N+t\Delta\uu)$ and $\dot\rr(t)=\sqrt{|\yy-\xx|}S(\yy,\xx,\Delta\uu)$. Noting that $g'$ is bounded and estimating as in \cite{BhattacharyaLipton} gives 
\begin{align*}
    |T(\xx,t)|\leq C||\uu^{N+1}-\uu^N||_\infty
\end{align*}
where $C$ will always denote a generic constant independent of $\uu^{N+1}-\uu^N$. From this we deduce
\begin{align}\label{equ: firstest}
    |\omega_{N+1}|\leq C||\uu^{N+1}-\uu^N||_\infty\leq||\vv^{N+1}-\vv^N||_\infty+||\UU^{N+1}-\UU^N||_\infty.
\end{align}


%
To conclude \eqref{equ: limiting1} we show $||\vv^{N+1}-\vv^N||_\infty\leq C ||\UU^{N+1}-\UU^N||_\infty$. Since one has that $|||\LL'_N[\uu^N]^{-1}||| \le \infty$ one sees from Banach's lemma and Lipschitz continuity \eqref{Lipschits} that there exists a closed ball $\overline{B(R,\uu^N)}=\{\uu:\, ||\uu-\uu^N||_\infty\leq R\}$  of radius $R$ and center $\uu^N$ for which $\LL'_N[\uu]^{-1}$ is well defined. and consequently, a positive constant $K_N>0$ determined at step $N$ for which $K_N||\vv||_\infty\leq ||\LL'[\uu]\vv||_\infty$ for any fixed $\uu$ in $\overline{B(R,\uu^N)}$ and for all $\vv\in\mathcal{V}$. We now can choose  $||\UU^{N+1}-\UU^N||_\infty$  sufficiently small so that $\uu_0=\vv^N+\UU^{N+1}$ lies inside $\uu$ in $\overline{B(R,\uu^0)}$. 
 
 Noting that $\mathcal{L}_{N}[\uu^N]=0$ and $\mathcal{L}_{N}[\uu^{N+1}]=0$ on ${D}_{N+1}$ gives
$$A=\mathcal{L}_N[\uu^{N+1}]-\mathcal{L}_N[\uu_0]=\mathcal{L}_N[\uu^N]-\mathcal{L}_N[\uu_0]=B$$ on $D_{N+1}.$ Here $\uu^{N+1}-\uu_0=\vv^{N+1}-\vv^N:=\Delta \vv^{N+1}$ and $\uu_0-\uu^N=\UU^{N+1}-\UU^N:=\Delta\UU^{N+1}$ and from the fundamental theorem of calculus as in \cite{BhattacharyaLipton},
\begin{align*}
  ||A||_\infty=  || \int_0^1{ \mathcal{L}}'_{N}[\boldsymbol{r}(t)]dt \Delta \vv ||_\infty
\end{align*}
with $\rr(t)=\sqrt{|\yy-\xx|}S(\yy,\xx,\uu^{N}+t\Delta\vv)$. From \eqref{lip}
\begin{align*}
    ||B||_\infty  \leq  C\,||\Delta \UU^{N} ||_\infty
\end{align*}
 Applying the mean value theorem and calculating as in \cite{Ortega} we get
\begin{align*}
  ||A||_\infty=  ||  \LL'_{N}[\boldsymbol{r}(\overline{t})]\,\Delta \vv ||_\infty,
\end{align*}
for some $0\leq \overline{t}\leq 1$. Since $\rr(\overline{t})$ is in the ball $\overline{B(R,\uu^N)}$ we get
\begin{align*}
 K_N||\Delta\vv|| _\infty\leq||A||_\infty,
\end{align*}
so
\begin{align*}
 K_N||\vv^{N+1}-\vv^N|| _\infty\leq||A||_\infty= ||B||_\infty  \leq  C\,||\UU^{N+1}-\UU^N ||_\infty
\end{align*}
where  $K_N$ is independent of $||\vv^{N+1}-\vv^N||_\infty$ and
 $||\UU^{N+1}-\UU^N||_\infty$ can be chosen small independently of $K_N$ so  \eqref{equ: limiting1} is established.

\section{Existence of nonlocal bond breaking evolution and emergence of a fracture evolution}
\label{exist}

We begin by stating two conditions that when taken together are sufficient  for the existence of an inverse of $\LL'[\uu]$ on  $\mathcal{V}$  by showing that $\LL'[\uu]$ satisfies the hypotheses of   (\cite{BhattacharyaLipton} Theorem 8).
One condition involves the stability tensor defined by: 
\begin{defn}{Stability tensor}\label{stabilitytensorload}
\begin{align}\label{stabtenload}
    \mathbb{A}[\uu](\xx) = \int\limits_{H_\epsilon(\xx) \cap \Omega}^{} \frac{J^\epsilon(\abs{\yy - \xx}) }{\epsilon^{n+1} \omega_n|\yy-\xx|}  g''\left( \sqrt{\abs{\yy - \xx}} S\left(\yy, \xx, \uu\right)\right) \ee_{\yy - \xx} \otimes \ee_{\yy - \xx}  d\yy.
\end{align}
\end{defn}
For a fixed bases $\mathbb{A}[\uu](\xx)$ is an $n\times n$ symmetric  matrix, $\mathbb{A}[\uu](\xx)=\mathbb{A}^T[\uu](\xx)$).
We write $\mathbb{A}^2[\uu]-\gamma^2\mathbb{I}>0$ when for all $\xx\in D$ and all $\vv\in\mathbb{R}^n$, $\mathbb{A}^2[\uu](\xx)\vv\cdot\vv-\gamma^2|\vv|^2>0$. The stability tensor for nonlocal modeling arises in many contexts  including dynamic fracture \cite{sillingweknerascaribobaru,lipton2014dynamic,lipton2019complex} and elasticity 
 \cite{DuGunLehZho,MengeshaDuNonlocal14,MengeshaDu}.  

To get the conditions, we set $\uu=\vv+\UU$ for $\vv$ in $\mathcal{V}$ and $\UU$ in $\mathcal{B}$ and integrate by parts as in (Lemma 5.3  of \cite{BhattacharyaLipton}) to find that $\LL'[\uu]$ is symmetric
on the $L^2(\Omega,\mathbb{R}^n)$ closure of $\mathcal{V}$. Moreover, from \cite{BhattacharyaLipton} we have that $\LL'[\uu]$ is a bounded operator on the $L^2(\Omega,\mathbf{R}^n)$ closure of  $\mathcal{V}$.
Collecting these results and applying Theorem 8 of \cite{BhattacharyaLipton} the sufficient conditions for an inverse are given by:
\\

\begin{propn}[Sufficient conditions for an inverse]
\label{inverse}

\noindent Given $\uu=\vv+\UU$ for boundary data $\UU$ in $\mathcal{B}$ and $\vv$ in  $\mathcal{V}$,  if $Ker\{\LL'[\uu]\}=\boldsymbol{0}$ for $\LL'[\uu]$ regarded as an operator on the $L^2(\Omega,\mathbb{R}^n)$ closure of $\mathcal{V}$ and if there exists $\gamma>0$ such that $\mathbb{A}^2[\uu]-\gamma^2\mathbb{I}>0$ then $\LL'[\uu]^{-1}$ is a bounded operator on $\mathcal{V}$.
\end{propn}
The condition that there exists $\gamma>0$ for which $\mathbb{A}^2[\uu]-\gamma^2\mathbb{I}>0$ is equivalent to saying that all eigenvalues of $\mathbb{A}[\uu]$ lie outside an interval about $0$.
In \cite{BhattacharyaLipton} it is shown that  $Ker\{\mathbb{A}[\uu](\xx)\}=\{\boldsymbol{0}\}$ on $D$ is a necessary condition for invertability.

Now we establish the existence of a bond breaking evolution described in Section \ref{bomdbreaking} from which a fracture evolution described in Section \ref{Fracture} can emerge.
This is done in two steps. First we prove that if the displacement $\uu^N$  is  a solution to the boundary value problem for load case $\UU^N$ and if it satisfies suitable hypotheses then one can find a solution $\uu^{N+1}$ to the boundary value problem for an appropriately chosen load case $\UU^{N+1}$; this is the statement of Proposition \ref{1stProp}.  In the following step we establish initialization. Here we show that the evolution starts with the initial choice $\uu^0=0$ and for this choice we  show it is possible to  apply Proposition \ref{1stProp}, to find a solution $\uu^1$ to the boundary value problem for an appropriately chosen load case $\UU^{1}$.
To conclude we give criteria for when the evolution terminates.\\

\noindent{\bf Load increment step}\\

We introduce a ball of radius $t^\ast$  surrounding a point $\vv^N$ in $\mathcal{V}$ given by $B(t^\ast,\vv^N)=\{\vv\in\mathcal{V}:\,||\vv-\vv^N||_\infty<t^\ast\}$ and denote its closure by  $\overline{B(t^\ast,\vv^N)}$. 
Now we state the existence theorem for a load increment in the evolution.
\begin{propn}[Solution for load step $N+1$]
\label{1stProp}
Given $\LL_N(\uu^N)=0$ for $\xx$ in $D$ and $\uu^N=\vv^N+\UU^N$ with $\vv^N$ in $\mathcal{V}$ and $\UU^N$ in $\mathcal{B}$.
If if there exists $\gamma>0$ such that $\mathbb{A}^2[\uu^N]>\gamma^2\mathbb{I}$ and $Ker\{\LL'[\uu^N]\}=\boldsymbol{0}$ for $\LL'[\uu]$ regarded as an operator on the $L^2(\Omega,\mathbb{R}^n)$ closure of $\mathcal{V}$ then there exists a load $\UU^{N+1}$ and initial deformation $\uu_0=\vv^N+\UU^{N+1}$ 
such that $\LL'[\uu_0]^{-1}$ exists and is a bounded linear transform on $\mathcal{V}$. Moreover there is a ball of radius $t^\ast$ centered at 
$\vv^N$ denoted by $B(t^\ast,\vv^N)$ such the unique fixed point  $\vv^{N+1}$ of the\\
Newton map $T(\vv):\overline{B(t^\ast,\vv^N)}\rightarrow\overline{B(t^\ast,\vv^N)}$
\begin{align}\label{newtonmap}
T(\vv)=\vv-\LL_N'[\vv+\UU^{N+1}]^{-1}\LL_N[\vv+\UU^{N+1}]
\end{align}
belongs to the  closed ball $\overline{B(t^\ast,\vv^N)}$. This fixed point is isolated as there exists a radius $t^{\ast\ast}>t^\ast$ and closed ball $\overline{B(t^{\ast\ast},\vv^N)}$ for which no other fixed point lies. Thus, $\LL_N[\uu^{N+1}]=0$ and breaking bonds according to \eqref{breakingbonds} with $\uu=\uu^{N+1}$ delivers the new operator $\LL_{N+1}[\vv+\UU]$ acting on $\mathcal{V}+\mathcal{B}$ and $\LL_{N+1}[\uu^{N+1}]=0$.
\end{propn}

We establish the proposition using the Newton-Kantorovich theorem \cite{Kantorovich,Ortega,GraggTapia}.  We have from \eqref{Lipschits} for $\UU^{N+1}$ fixed and $\vv$ in $\mathcal{V}$
\begin{align}\label{Lipschitzz}
    |||\LL'_N[\vv+\UU^{N+1}]-\LL'_N[\uu_0]||| &<  K||\vv-\vv^N||_\infty.
\end{align}
From the hypothesis of Proposition \ref{1stProp} and Proposition \ref{inverse}  we have a finite $\beta$ such that
\begin{align}\label{normofinverse}
    |||\LL'_N[\uu_0]^{-1}]||| &<  \beta.
\end{align}
Now $\LL_N[\uu^N]=0$ and from \eqref{lip}
\begin{align}\label{prodcut}
    ||\LL_N[\uu^N]-\LL_N[\uu_0]||_\infty\leq C ||\UU^{N+1}-\UU^N||_\infty,
\end{align}
so we can choose $\UU^{N+1}$  so that $||\UU^{N+1}-\UU^N||_\infty$ is sufficiently small so that $||\LL_N'[\uu_0]^{-1}\LL_N[[\uu_0||_\infty \leq\eta$ with $h=\beta K \eta\leq 1/2$.
Then from the Newton-Kantorovich theorem we can choose radii  $t^\ast$ and $t^{\ast\ast}$ such that
\begin{equation}\label{radii}
t^\ast=\frac{1}{\beta K}(1-\sqrt{1-2h}),\qquad t^{\ast\ast}=\frac{1}{\beta K}(1+\sqrt{1-2h}),
\end{equation}
and the Newton iterates converge to the unique fixed point $\vv^{N+1}$ in the closed ball  $\overline{B(t^\ast,\vv^N)}$. Moreover, $\LL_N[\uu^{N+1}]=0 $ and breaking bonds  according to \eqref{breakingbonds} with $\uu=\uu^{N+1}$ delivers the new operator $\LL_{N+1}[\vv+\UU]$ acting on $\mathcal{V}+\mathcal{B}$ and $\LL_{N+1}[\uu^{N+1}]=0$. Last, the fixed point is isolated and no other fixed points lie inside the larger closed ball  $\overline{B(t^{\ast\ast},\vv^N)}$.\\

\noindent{\bf Initial load step}\\

We now start with the initial choice $\uu^0=0$ and show that there exists a solution $\uu^1$ to the boundary value problem for appropriate nonzero boundary data $\UU^1$.  

\begin{propn}
\label{initialize}
On choosing  $\uu^0=0$ we have that $\LL'[\uu^0]^{-1}$ exists and is a bounded linear operator mapping $\mathcal{V}$ onto $\mathcal{V}$ moreover there is a radius $R>0$
for which  $\uu^0=0$ is the unique solution of the boundary value problem
\begin{equation} 
\LL[\uu^0](\xx)=0, \hbox{ for $\xx$ in $D$ and $\uu^0=0$ on $\Omega_d$ },
\end{equation}
among all functions in the ball $B(R,0)$.
\end{propn}
 Proposition \ref{initialize} shows that the initialization $\uu^0=0$, satisfies the hypotheses of Proposition \ref{1stProp} so there is a least one  choice of nonzero  boundary load  $\UU^1$  for which there is a solution $\uu^1$ of the boundary value problem \eqref{eq:bvp}. Breaking bonds if needed according to Proposition \ref{1stProp} gives the new operator $\LL_1$ on $\mathcal{V}+\mathcal{B}$ and $\LL_1[\uu^1](\xx)=0$ for $\xx$ in $D$. The proof of Proposition \ref{initialize}  is given in the Appendix.
\\

\noindent{\bf Terminal load step}\\

The load increment is terminated when $Ker\{\LL'_{N+1}[\uu^{N+1}]\}\not=0$ or when $\mathbb{A}[\uu^{N+1}]$ vanishes on a set of finite volume.
Here it is found that $\LL'_{N+1}[\uu^{N+1}]^{-1}$ does not exists when the stability tensor  $\mathbb{A}[\uu^{N+1}](\xx)$ vanishes on a set of finite volume  inside $D$, see \cite{BhattacharyaLipton}.


%

\section{Numerical algorithm for fracture evolution}
\label{Newt}
The solution $\uu^N$ to $\LL_N[\uu^N]=0$ for each loading $\UU^N$  in the fracture evolution is found numerically using the Newton-Raphson method. We prescribe an increment of boundary load $\Delta\UU=\UU^{N}-\UU^{N-1}$ and
starting with the initial guess $\uu^N_0=\vv^{N-1}+\UU^{N}$, for all $k$ solve for $\Delta \uu$
\begin{align}
\label{eq:newton}
    -\LL'_{N-1}[\uu^N_k] \Delta \uu = \LL_{N-1}[\uu^N_k] ,
\end{align}
and set $\uu^N_{k+1} = \uu^N_k + \Delta \uu$. The approximate root $\uu^{N}$ satisfies
\begin{align}
\LL_{N-1}[\uu^{N}](\xx)\approx 0, 
\end{align}
and $\uu^{N}=\vv^{N}+\UU^{N}=\UU^{N}$ on $\Omega_d$.
In the numerical implementation, we approximate and break all bonds according to 
\begin{align}\label{eq: brokenbusted}
    {r^e}>{r=\sqrt{|\yy - \xx|}}S(\yy, \xx, \uu) \hbox{ or }
  \sqrt{|\yy-\xx|} S(\yy, \xx, \uu)&= r > {r^c}.
\end{align}
and $\LL_N[\uu^N](\xx)\approx 0$. 
Here the approximate equality ``$\approx$''  is measured by $||\LL_{N-1}[\uu^N]||_\infty\leq \tau$ 
with $\tau$ a preset tolerance.
The initial guess for the next increment $\UU^{N+1}$ is $\uu_0=\vv^{N}+\UU^{N+1}$ and the process is repeated.

\subsection{Discretization}
\label{sub:discretization}
Expanding $S(\yy, \xx, \ww) = \frac{\ww(\yy) - \ww(\xx)}{\abs{\yy - \xx}} \cdot \frac{\yy - \xx}{\abs{\yy - \xx}}$ we can write
\begin{align}
    \label{eq:expand}
 &   (L_{N-1}'[\uu] \ww (\xx))_i = \nonumber\\
 &   \sum_{j=1}^2 
    \left(
    \int\limits_{H_\epsilon(\xx) \cap \Omega_{N-1}}^{} \sigma_{ij}(\yy , \xx, \uu) w_j(\yy) d\yy - w_j(\xx) \int\limits_{H_\epsilon(\xx) \cap \Omega_{N-1}}^{} \sigma_{ij}(\yy , \xx, \uu) d\yy
    \right)
\end{align}	
where
\begin{align*}
    \sigma_{ij}(\yy, \xx, \uu)  = \frac{J^\epsilon\left( \abs{\yy - \xx}  \right) }{\epsilon^{n+1} \omega_n \abs{\yy - \xx}} g''(\sqrt{\yy - \xx} S(\yy, \xx, \uu)) \frac{y_i - x_i}{\abs{\yy - \xx}} \frac{y_j - x_j}{\abs{\yy - \xx}}.
\end{align*}	
The domain $\Omega$ is discretized at finitely many points. Let the set of points within the interior $D_{N-1}$ be given by $\{ \bx_i \}_{i=1}^M$.  Let $V_i$ be the volume element associated with the point $\bx_i$. 
For $k = 1, \dots, M$, denote the set of neighboring indices by $I_k^{N-1} = \{ l : \abs{\bx_k - \bx_l} \le \epsilon,\ \bx_l \in \Omega_{N-1},\ l \ne k \}$.
Define $W_i^k : = w_i(\bx_k)$ and $U_i^k : = u_i(\bx_k)$.

The discrete version of \eqref{eq:expand} is given by the linear operator $\K$ acting on the $2M$-dimensional vector $\mathbf{W} =\{ W_i^k : i =1,2 \text{ and } k=1, \dots, M  \}$
as
\begin{align*}
    (\K \mathbf{W})_i^k = 
    \sum_{j=1}^{2} \sum_{l \in I_k^{N-1}}^{} \sigma_{ij}(\bx_l, \bx_k, \bu) W_j^l V_l - \sum_{j=1}^{2} W_j^k \sum_{l \in I_k^{N-1}}^{} \sigma_{ij}(\bx_l, \bx_k, \bu) V_l,
\end{align*}	
where $\mathbf{U} =\{ U_i^k : i =1,2 \text{ and } k=1, \dots, M  \}$.

Similarly, denoting $B_i^k := (\LL_{N-1}[\uu](\bx_k))_i$ the discrete description of $\LL_{N-1}[\uu]$ is given by
\begin{align*}
    B_i^k & = \frac{2}{\epsilon^{n+1} \omega_n} 
    \sum_{l \in I_k^{N-1}}
    \frac{J^\epsilon(\abs{\bx_l - \bx_k})}{\sqrt{\abs{\bx_l - \bx_k}} } 
    g' \left( \sqrt{\abs{\bx_l - \bx_k}} S\left(\bx_l, \bx_k, \bu \right)  \right)  \frac{X_{l,i} - X_{k,i}}{\abs{\bx_l - \bx_k}} V_l.
\end{align*}	
Solving the $2M \times 2M$ system of linear equations given by
\begin{align*}
    -\K \mathbf{W} = \mathbf{B}
\end{align*}	
with $\uu = \uu_k^{N-1}$
one obtains the discretization $\mathbf{W}$ of the solution $\Delta \uu$ to \eqref{eq:newton}.
The Newton iteration is terminated with the residual $\norm{\LL[\uu^{N-1}_k]}_\infty$ is below a prescribed tolerance.

To prevent overshooting and to improve the stability of the Newton's iteration, a successive over-relaxation is employed where the solution increment $\uu_{k+1}^N = \uu_k^N + \Delta \uu$ is modified to $\uu_{k+1}^N = \uu_k^N + \theta \Delta \uu$, with $\theta = 0.5$.

\section{Numerical examples}
\label{sec:numerics}

Here, we consider a set of numerical examples and recover stable crack growths under displacement controlled loading.
We consider a straight crack growth under uniaxial tension. We observe a reduction of energy in the intact material once the cracks starts to grow.
Inspired by experiment, we consider an L-shaped panel subjected to an in-plane torsional loading. As expected, the crack starts to grow from the re-entrant corner.
Finally, we consider a rectangular panel with two pre-notches with variable offsets. Under uniaxial tension, the cracks start to grow inward. 
The offset distance of the pre-cracks affects the interaction of the growing cracks.

In the discrete setting, the quantity \textit{damage} is defined as
\begin{align*}
    d(\xx) = 1 - \frac{\#\text{intact bonds connected to $\xx$}}{\#\text{total bonds connected to $\xx$ in the reference configuration}}.
\end{align*}
We give two types of plots for the fracture evolution. 
Figures \ref{fig:4a} and \ref{fig:7a} display intact bonds (straight line segments) and in Figures \ref{fig:4b}, \ref{fig:7b}, and \ref{fig:double-notched} we plot the damage $d$.

\subsection{Straight crack propagation}
\label{sec:straight}

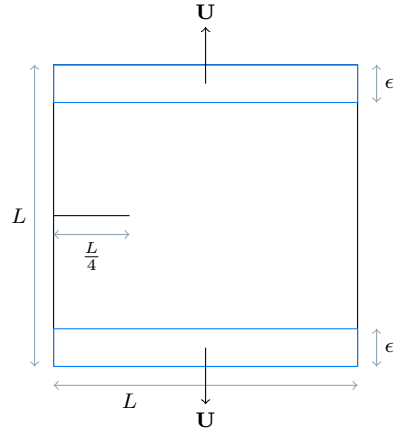
\begin{figure}[tb]
    \centering
    \begin{tikzpicture}
    \draw (0,0) rectangle (4,4);
    \draw (0,2) -- (1,2);
    \draw[azure] (0,0) rectangle (4,0.5);
    \draw[azure] (0,3.5) rectangle (4,4);
    \draw[<->, cadetgrey] (0,-0.25) -- (4,-0.25);
    \node[below] at (1,-0.25) {$L$};
    \draw[<->, cadetgrey] (-0.25,0) -- (-0.25,4);
    \node[left] at (-0.25, 2) {$L$};
    
    \draw[<->, cadetgrey] (0,1.75) -- (1,1.75);
    \node[below] at (0.5, 1.75) {$\frac{L}{4}$};
    
    \draw[<->, cadetgrey] (4.25,3.5) -- (4.25,4);
    \node[right] at (4.25, 0.25) {$\epsilon$};
    \draw[<->, cadetgrey] (4.25,0) -- (4.25,0.5);
    \node[right] at (4.25, 3.75) {$\epsilon$};
    
    \draw[->] (2,0.25) -- (2,-0.5);
    \node[below] at (2, -0.5) {$\UU$};
    \draw[->] (2,3.75) -- (2,4.5);
    \node[above] at (2, 4.5) {$\UU$};
    \end{tikzpicture}
    \caption{Square domain with a single horizontal pre-notch}
    \label{fig:straight-diagram}
\end{figure}

\begin{figure}
    \centering
    \subfloat[Bonds]{
\includegraphics[width=0.3\linewidth]{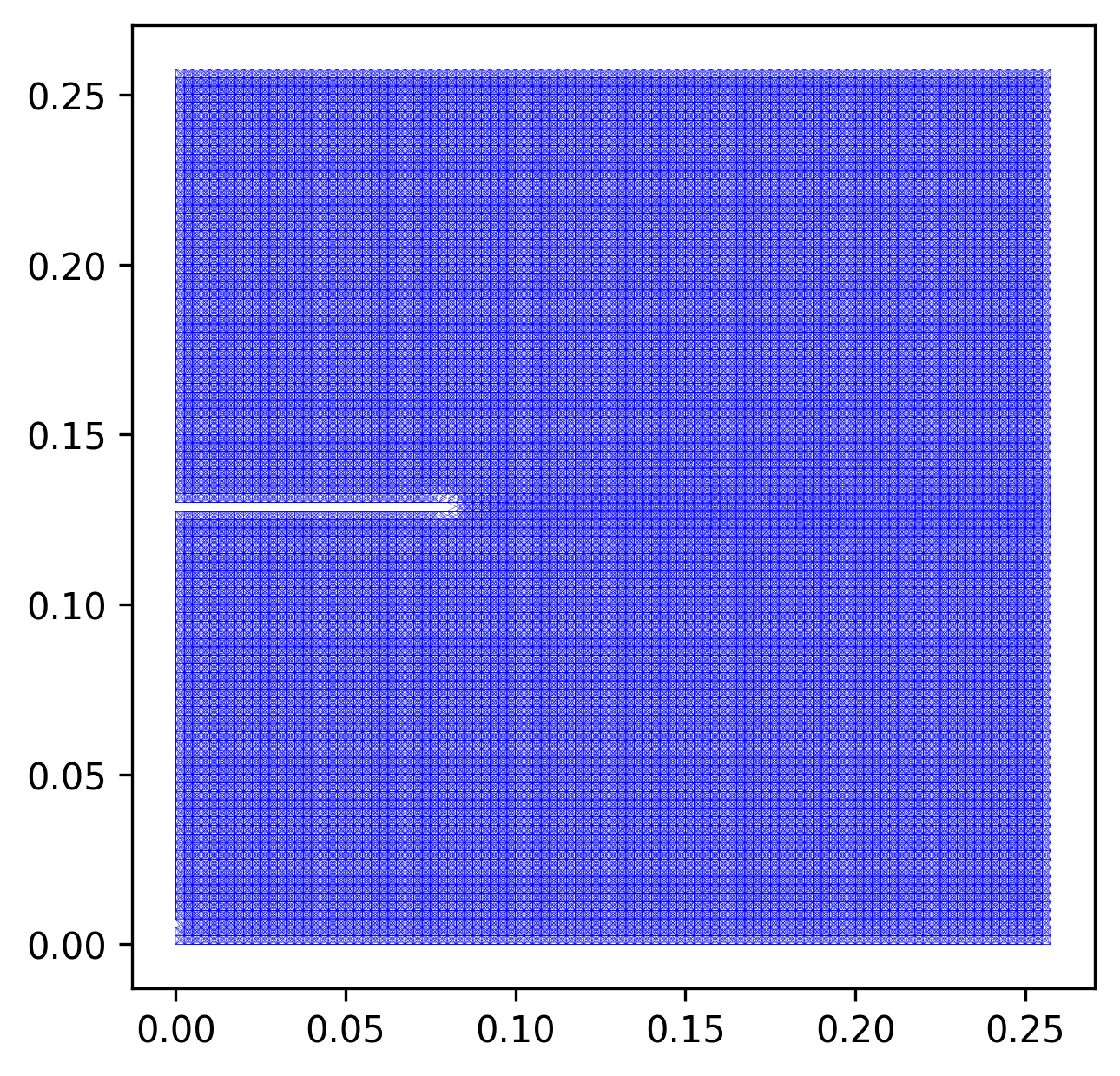}
\includegraphics[width=0.3\linewidth]{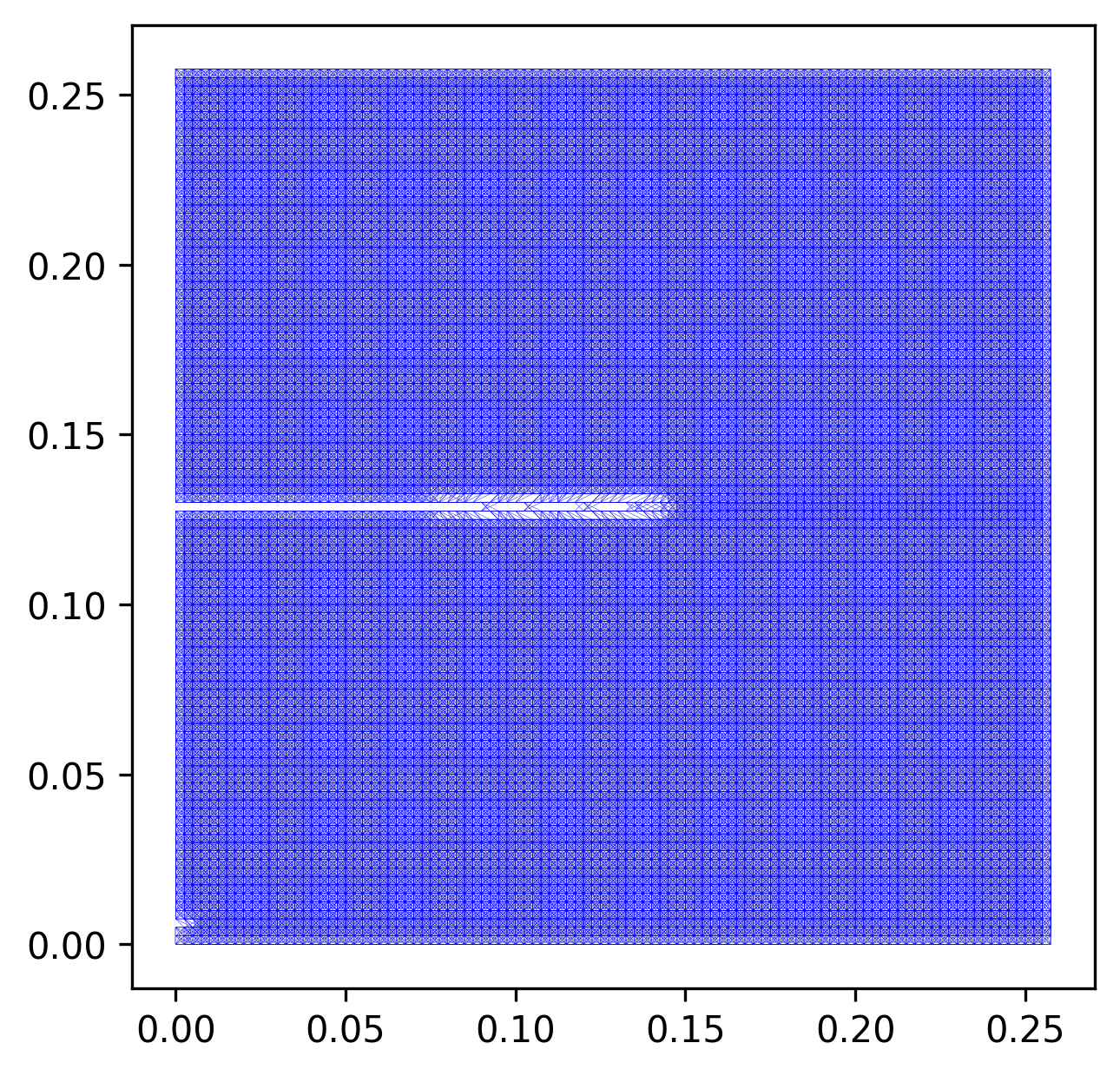}
\includegraphics[width=0.3\linewidth]{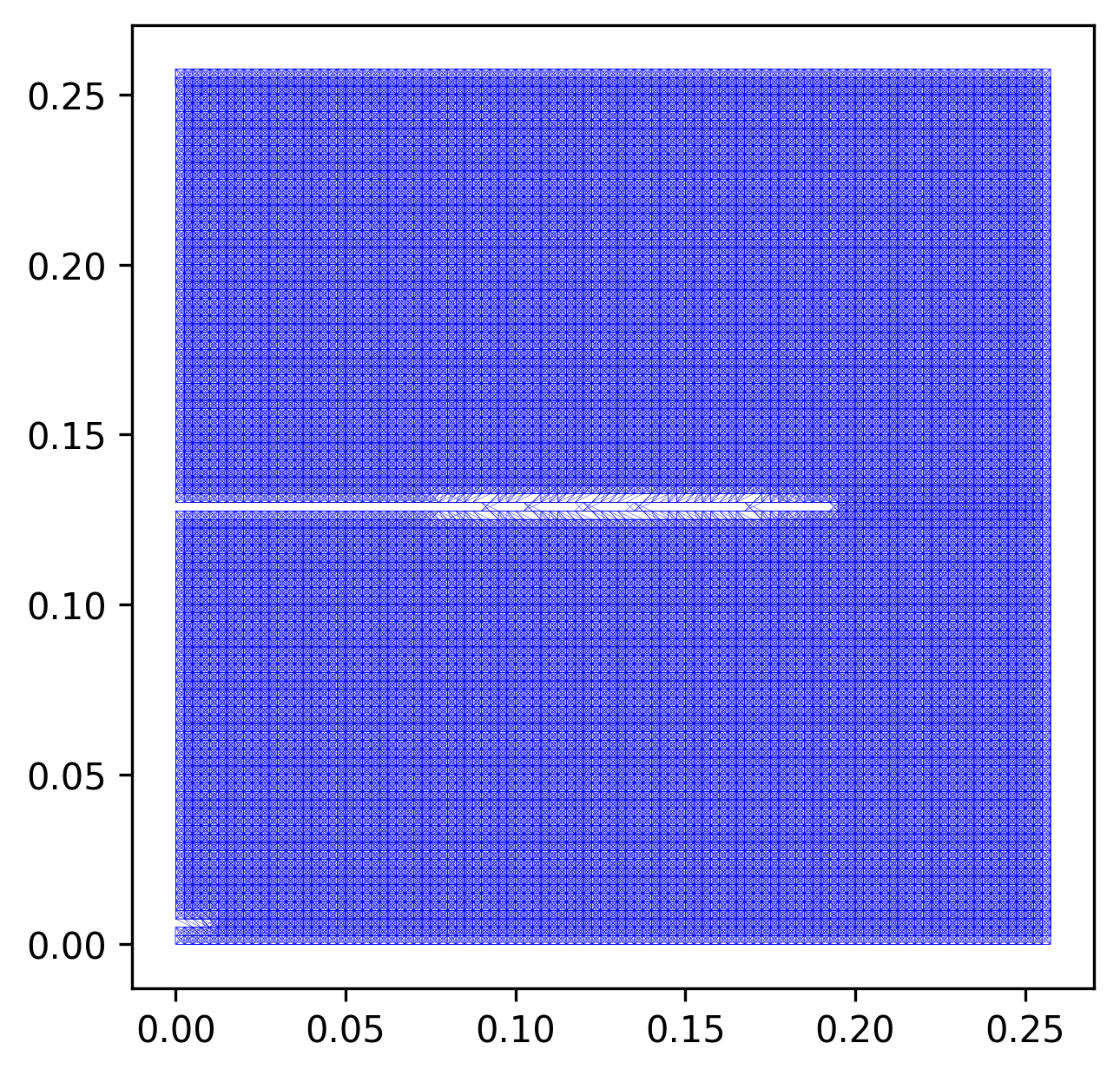}
\label{fig:4a}
}
\\
\subfloat[Damage]{
\includegraphics[width=0.3\linewidth]{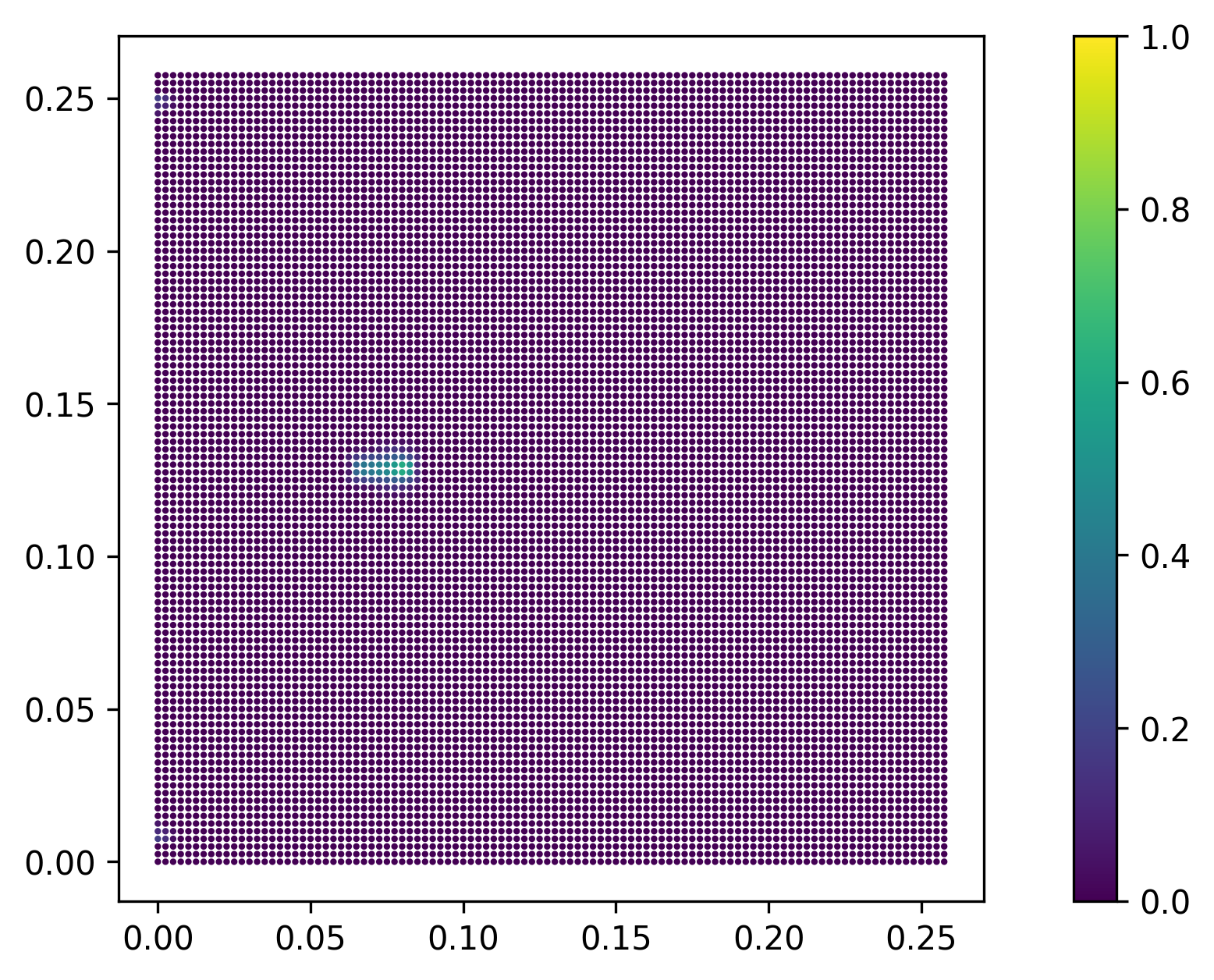}
\includegraphics[width=0.3\linewidth]{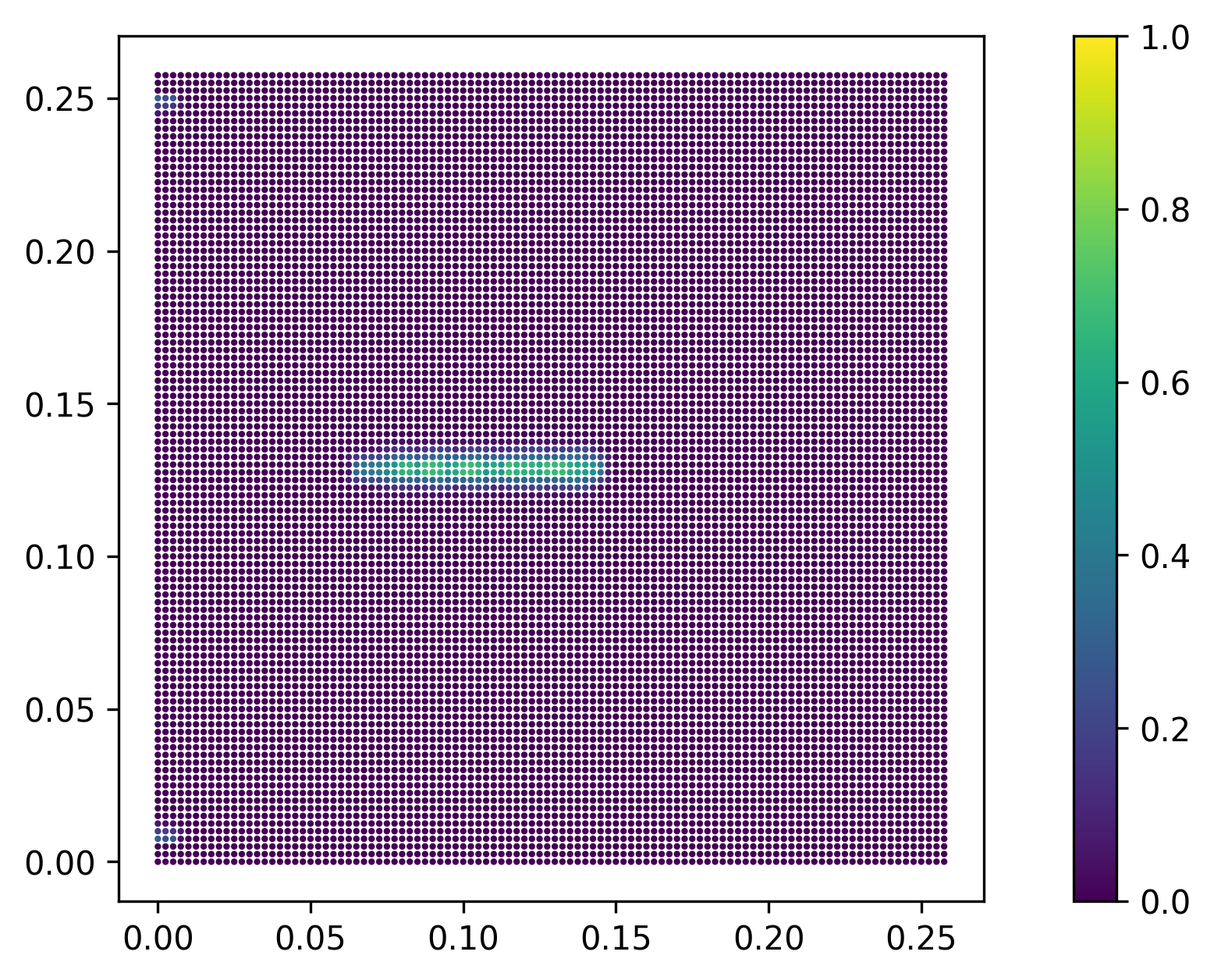}
\includegraphics[width=0.3\linewidth]{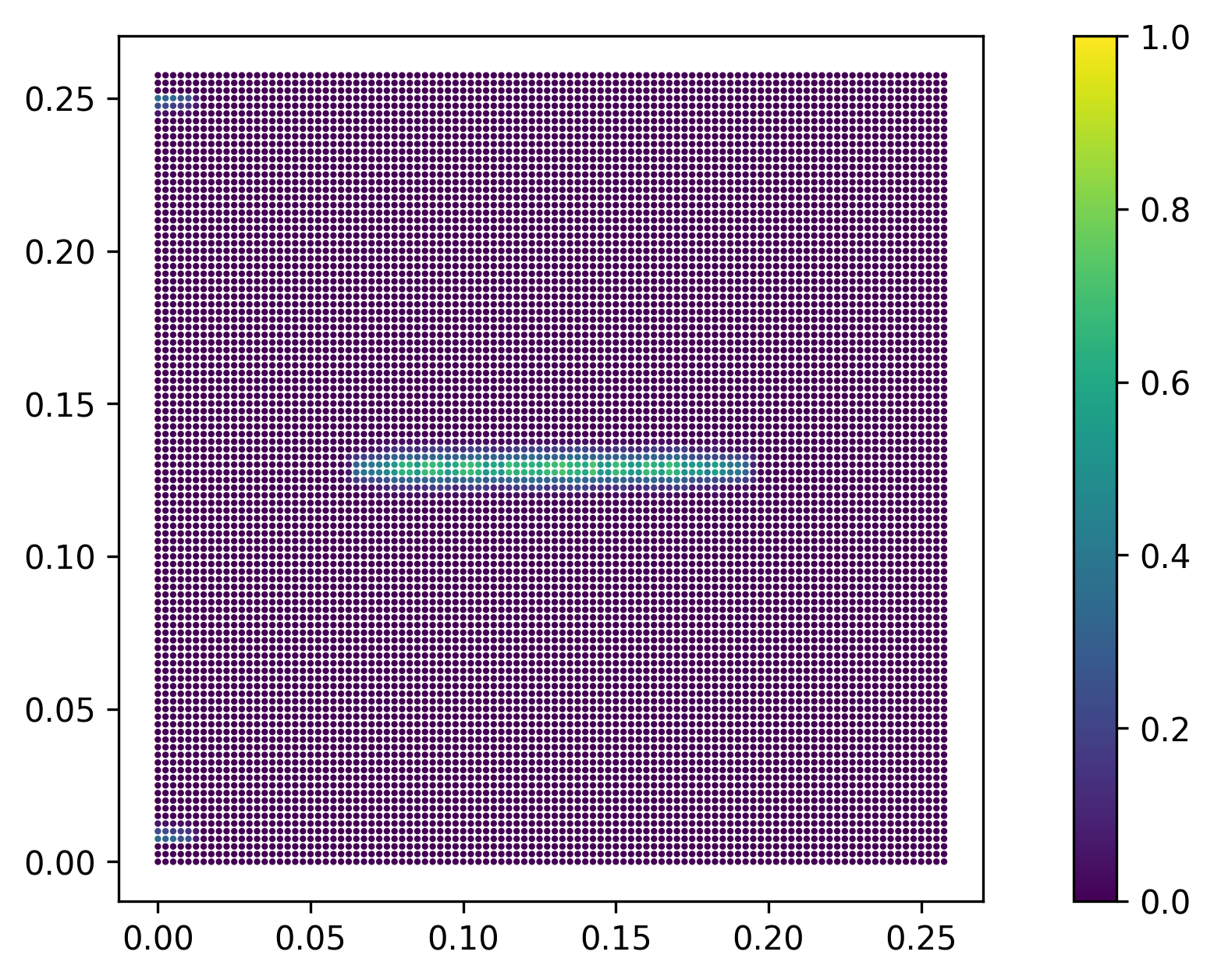}
\label{fig:4b}
}
\\
\subfloat[Energy]{
\label{fig:straight-energy}
\includegraphics[width=0.6\linewidth]{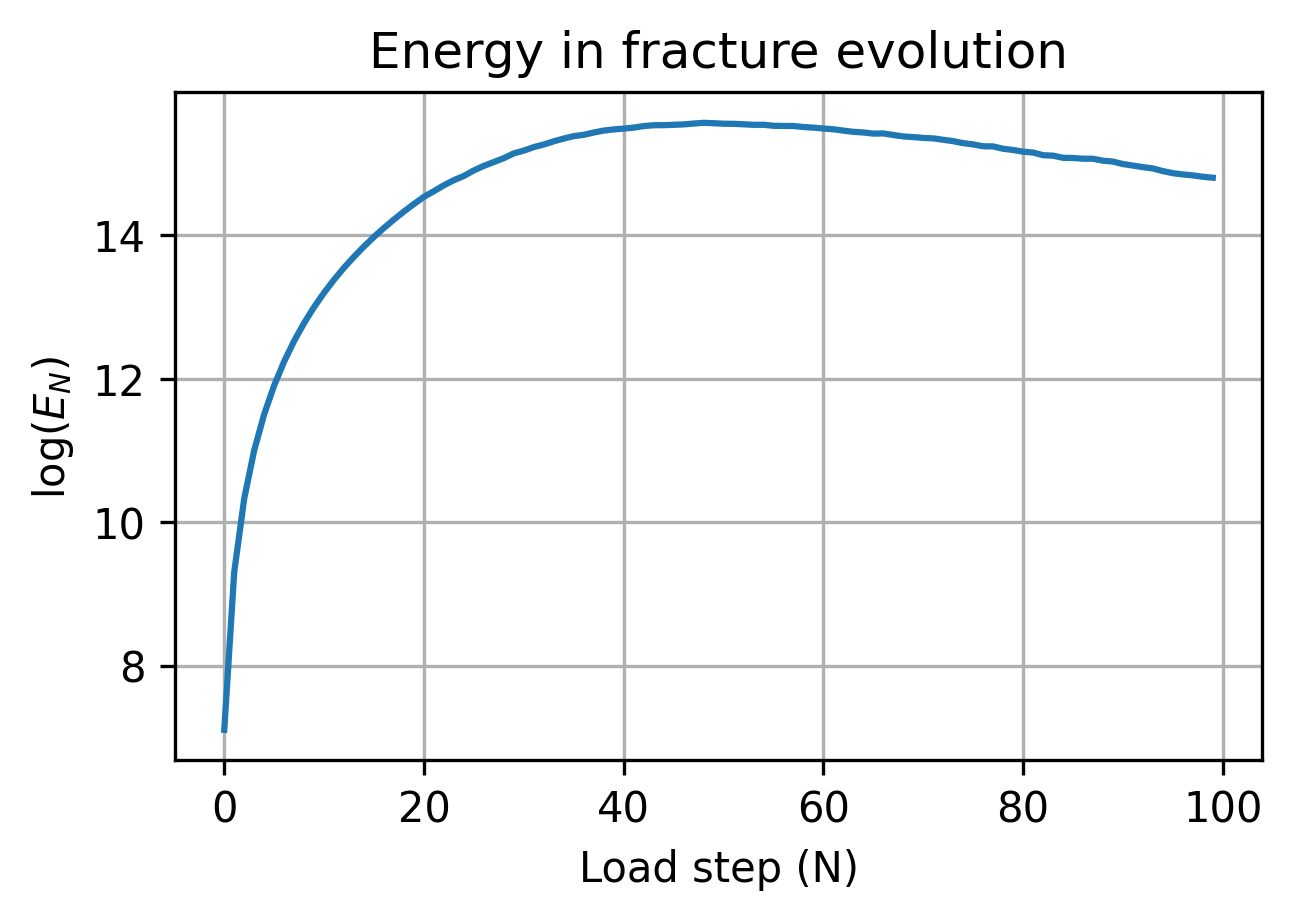}
}
\\
\subfloat[Residual]{
\label{fig:straight-residue}
\includegraphics[width=0.6\linewidth]{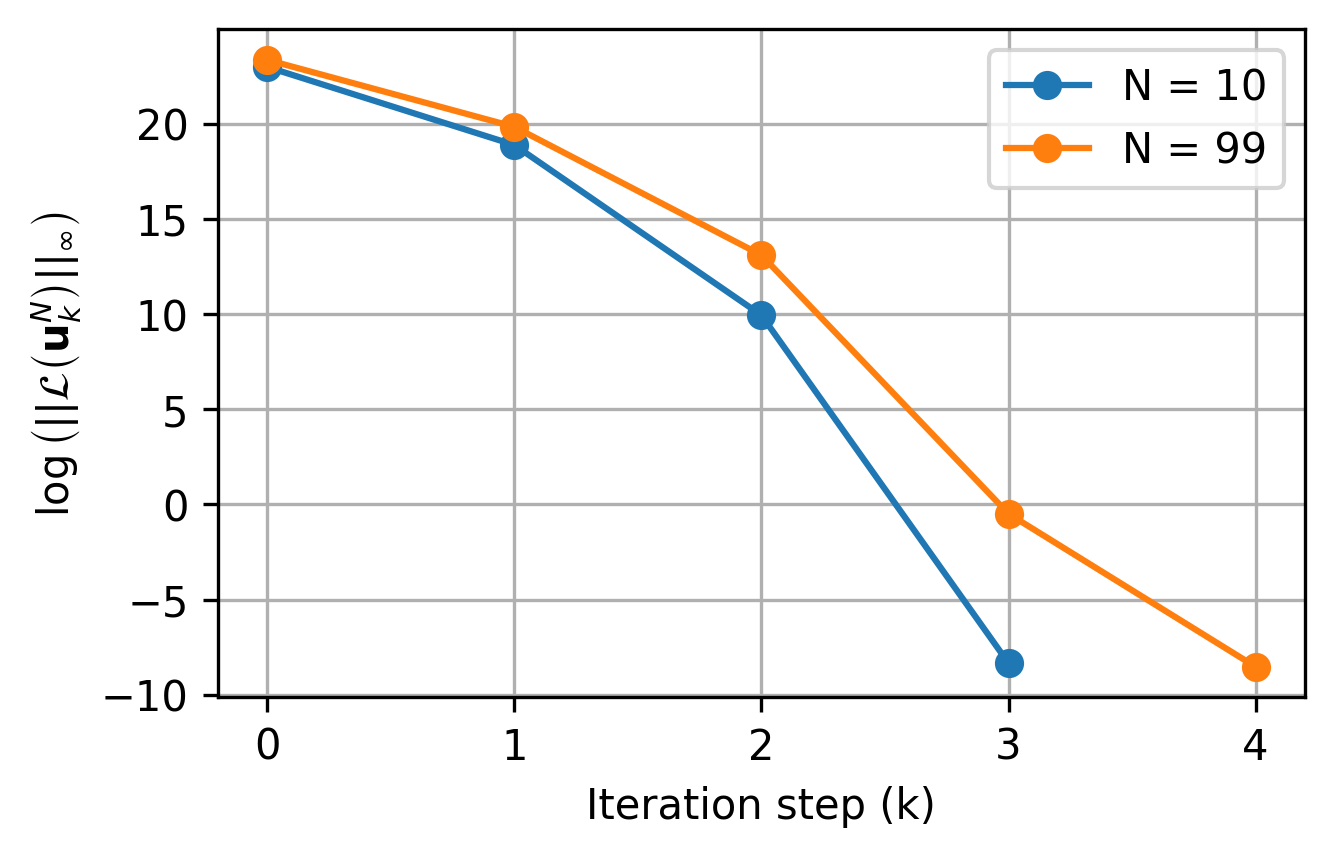}
}
		  \caption{Straight crack propagation in a square domain}
		  \label{fig:straight}
\end{figure}

We consider a square-shaped domain length $L$ with a horizontal pre-notch of length $\frac{L}{4}$ on the left edge of the domain, where L = 260 mm (see Fig. \ref{fig:straight-diagram}). A regular rectangular grid is considered with meshsize $h = 2.5$ mm. The peridynamic horizon is taken to be $\epsilon =$ 3h.
The material parameters are taken to be $E = 210$ GPa and $G_c = $ 2700 J/m$^2$.
The domain is subjected to uniaxial tension in the vertical direction via a displacement on the nonlocal boundary, which is taken to be a layer of thickness $\epsilon$ adjacent to the top and the bottom edge of the rectangle. We apply a vertical outward displacement of $U^0 = 2.05 \times 10^{-5}$ mm.
The displacement is incremented by $\Delta U = 4.12 \times 10^{-5}$ mm at each load step.
The crack grows at the tip of the pre-notch and extends horizontally in a straight line.
Fig. \ref{fig:straight} shows the crack path evolution at load steps $N=40, 70$, and $99$.
The energy $E_N$ associated with the quasistatic evolution is shown in Fig. \ref{fig:straight-energy}.
In Fig. \ref{fig:straight-residue}, we plot the residuals $\norm{\LL[\uu^N_k]}_{L^\infty(\Omega_d)}$ associated with the Newton iteration for two different load steps. We plot the residual versus Newton step for a load step before the crack starts (N = 10) and for a load step after the crack has progressed significantly (N = 99). The residuals lie below a tolerance of $10^{-5}$ after $3$ iterations for load step $(N=10)$ and after $4$ iterations for load step $(N=99)$.  

\subsection{L-shaped panel test}
\label{sec:L-shaped}

\begin{figure}
    \centering
    \includegraphics[width=0.6\linewidth]{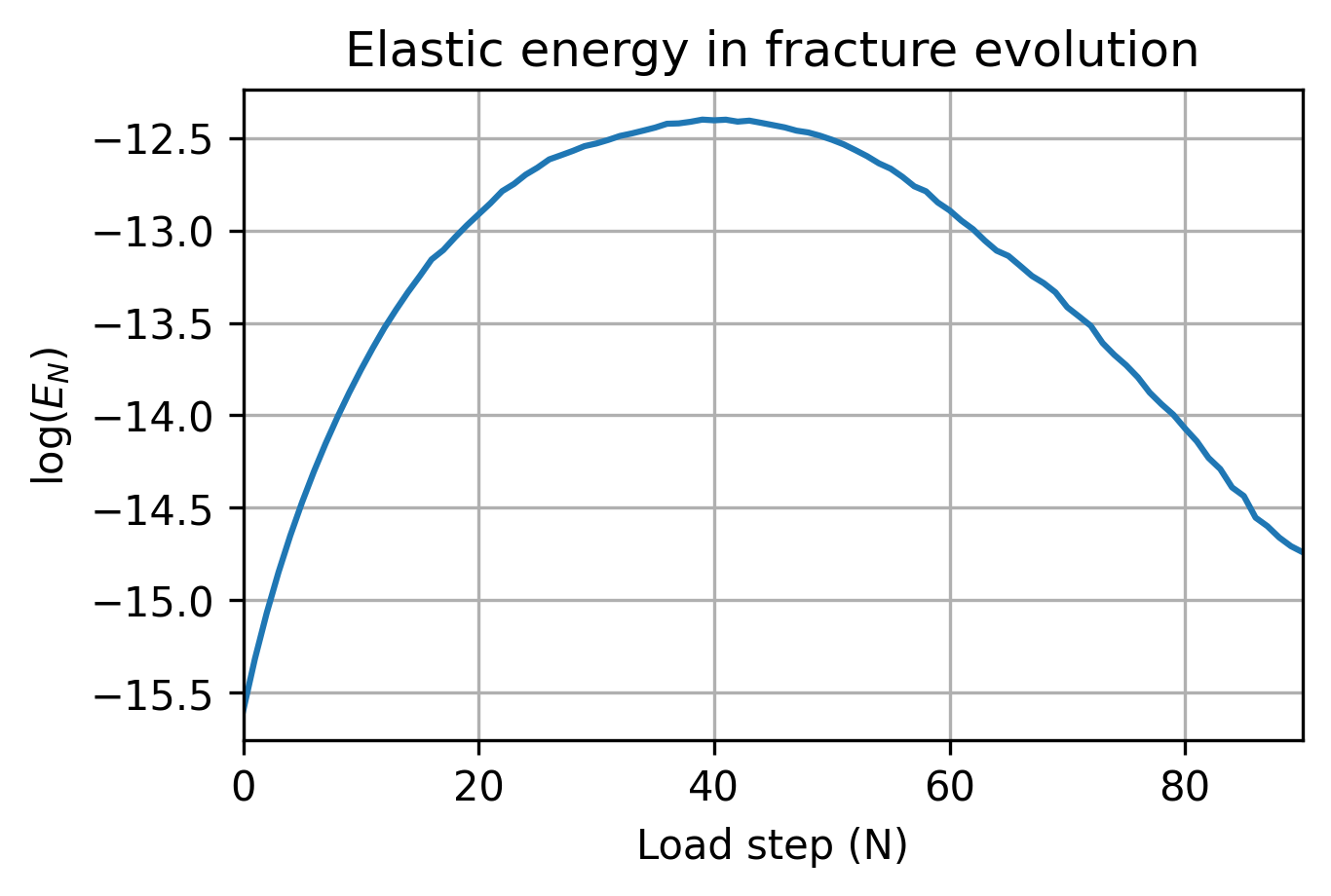}
		  \caption{Energy versus the load step for the L-shaped panel test.}
		  \label{fig:energy}
\end{figure}

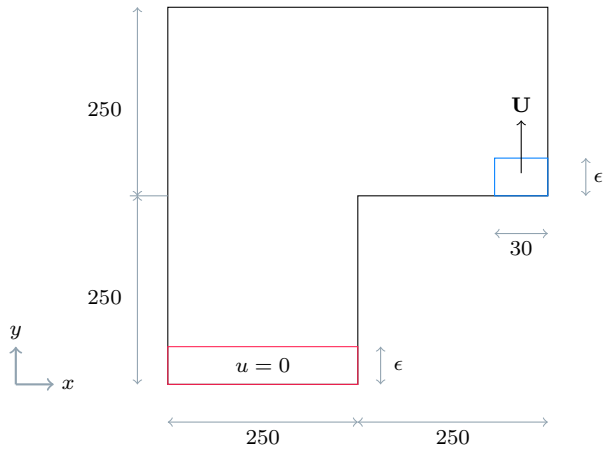
\begin{figure}[tb]
    \centering
    \begin{tikzpicture}
    \draw (0,0) -- (2.5,0) -- (2.5,2.5) -- (5,2.5) -- (5,5) -- (0,5) -- cycle;
    \draw[awesome] (0,0) rectangle (2.5,0.5);
    \node at (1.25, 0.25) {$u = 0$};
    \draw[azure] (4.3, 2.5) rectangle (5, 3);
    \draw[<->, cadetgrey] (4.3, 2.0) -- (5, 2.0);
    \node[below] at (4.65, 2.0) {$30$};
    \draw[<->, cadetgrey] (5.5, 2.5) -- (5.5, 3.0);
    \node[right] at (5.5, 2.75) {$\epsilon$};
    \draw[->]  (4.65,2.8) -- (4.65,3.5);
    \node[above] at (4.65, 3.5) {$\UU$};
\draw[<->,cadetgrey]  (2.8,0) -- (2.8,0.5);
\node[left] at (3.2,0.25) {$\epsilon$};
%
%
\draw[cadetgrey] (-0.5,2.5) -- (0,2.5);
\draw[<->,cadetgrey]  (-0.4,0) -- (-0.4,2.5);
\node[left] at (-0.5,1.15) {$250$};
\draw[<->,cadetgrey]  (-0.4,2.5) -- (-0.4,5);
\node[left] at (-0.5,3.65) {$250$};
\draw[<->,cadetgrey]  (0,-0.5) -- (2.5,-0.5);
\node[below] at (1.25,-0.5) {$250$};
\draw[<->,cadetgrey]  (2.5,-0.5) -- (5,-0.5);
\node[below] at (3.75,-0.5) {$250$};

\draw[->,thick,cadetgrey] (-2,0) -- (-1.5,0);
\draw[<-,thick,cadetgrey] (-2,0.5) -- (-2,0);
\node[above] at (-2,0.5) {\small $y$};
\node[right] at (-1.5,0) {\small $x$};

    \end{tikzpicture}
    \caption{L-shaped domain geometry and loading}
    \label{fig:diagram-L}
\end{figure}

\begin{figure}
    \centering
    \subfloat[Bonds]{
\includegraphics[width=0.3\linewidth]{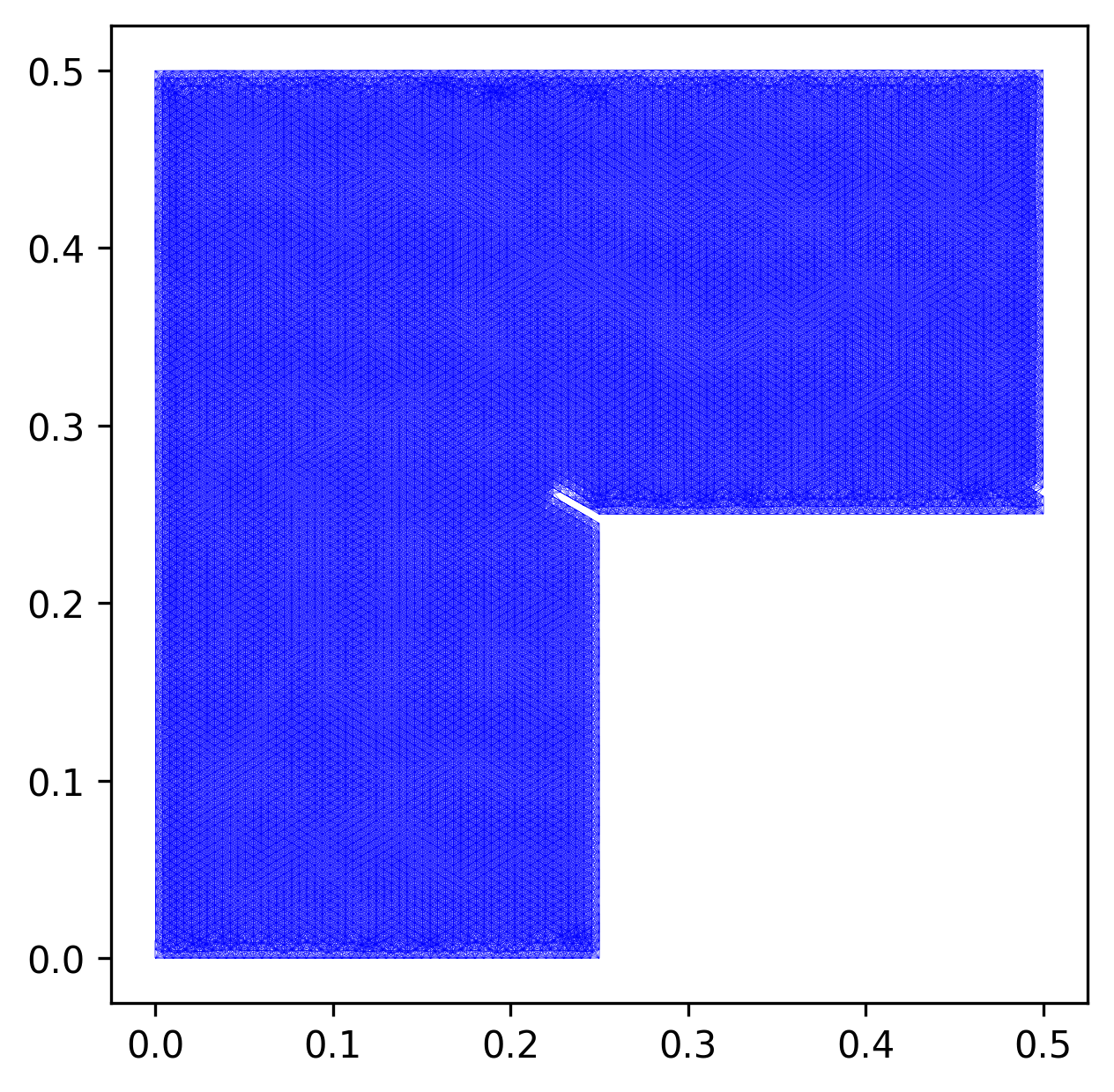}
\includegraphics[width=0.3\linewidth]{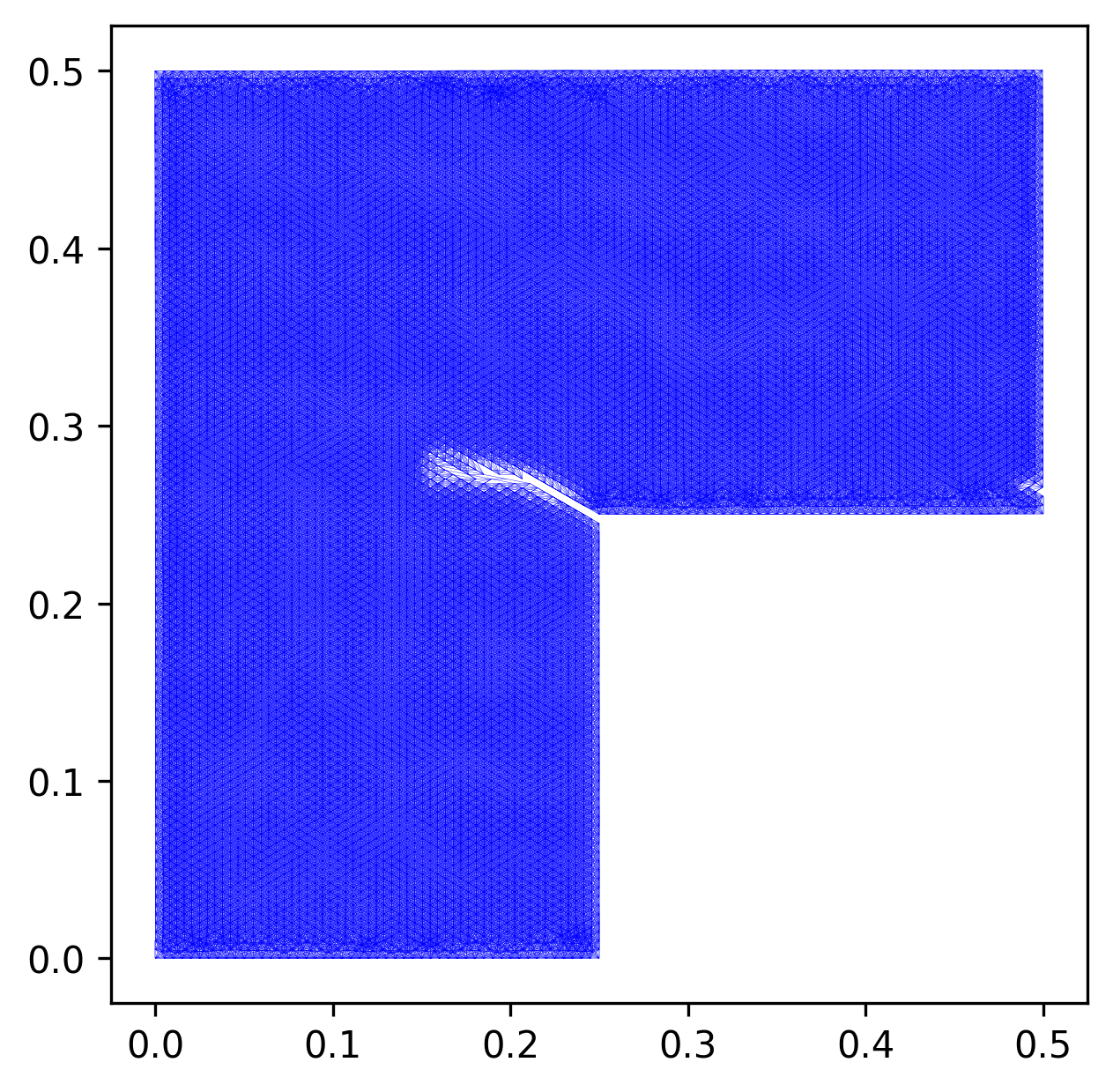}
\includegraphics[width=0.3\linewidth]{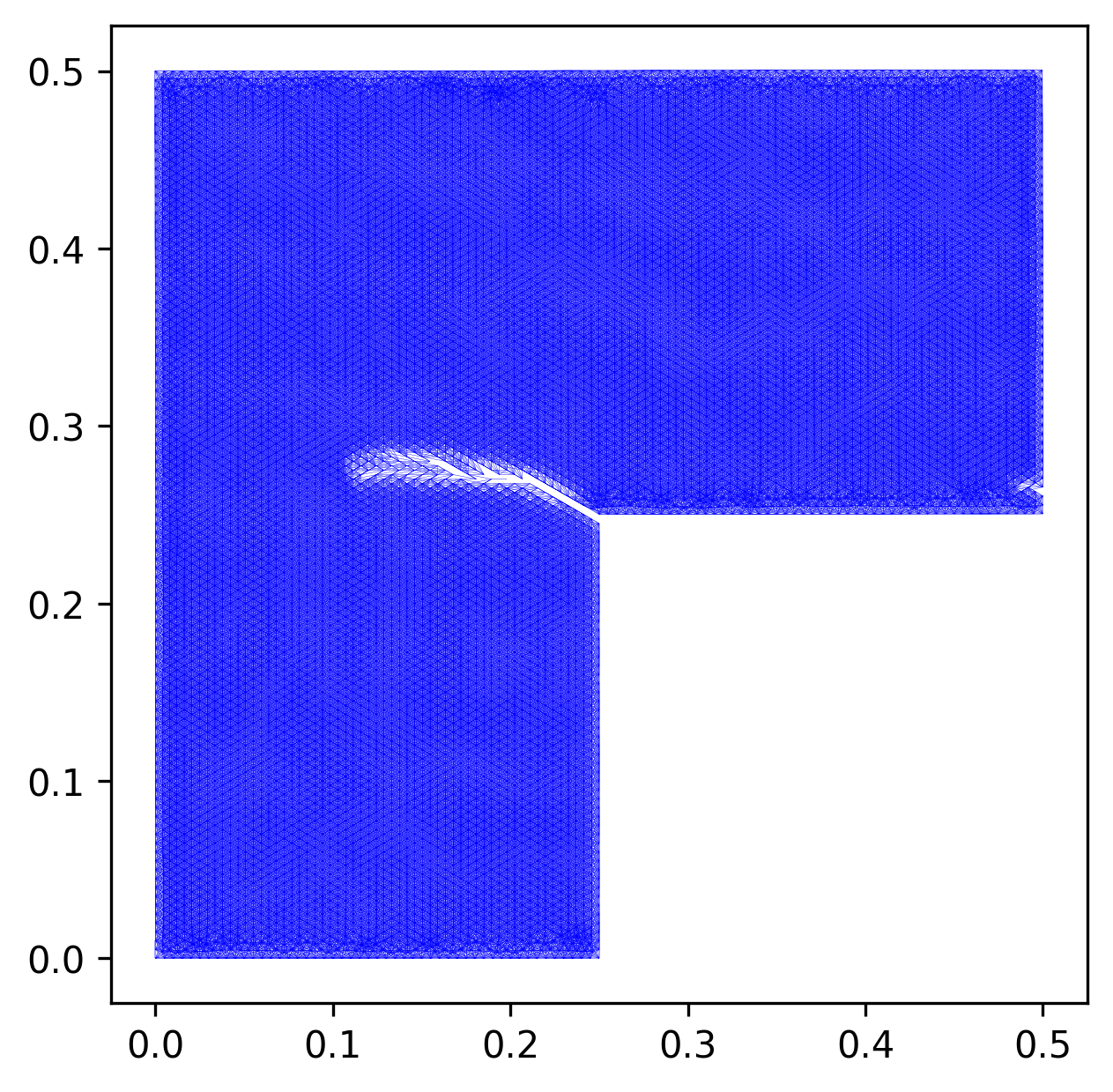}
\label{fig:7a}
}
\\
\subfloat[Damage]{
\includegraphics[width=0.3\linewidth]{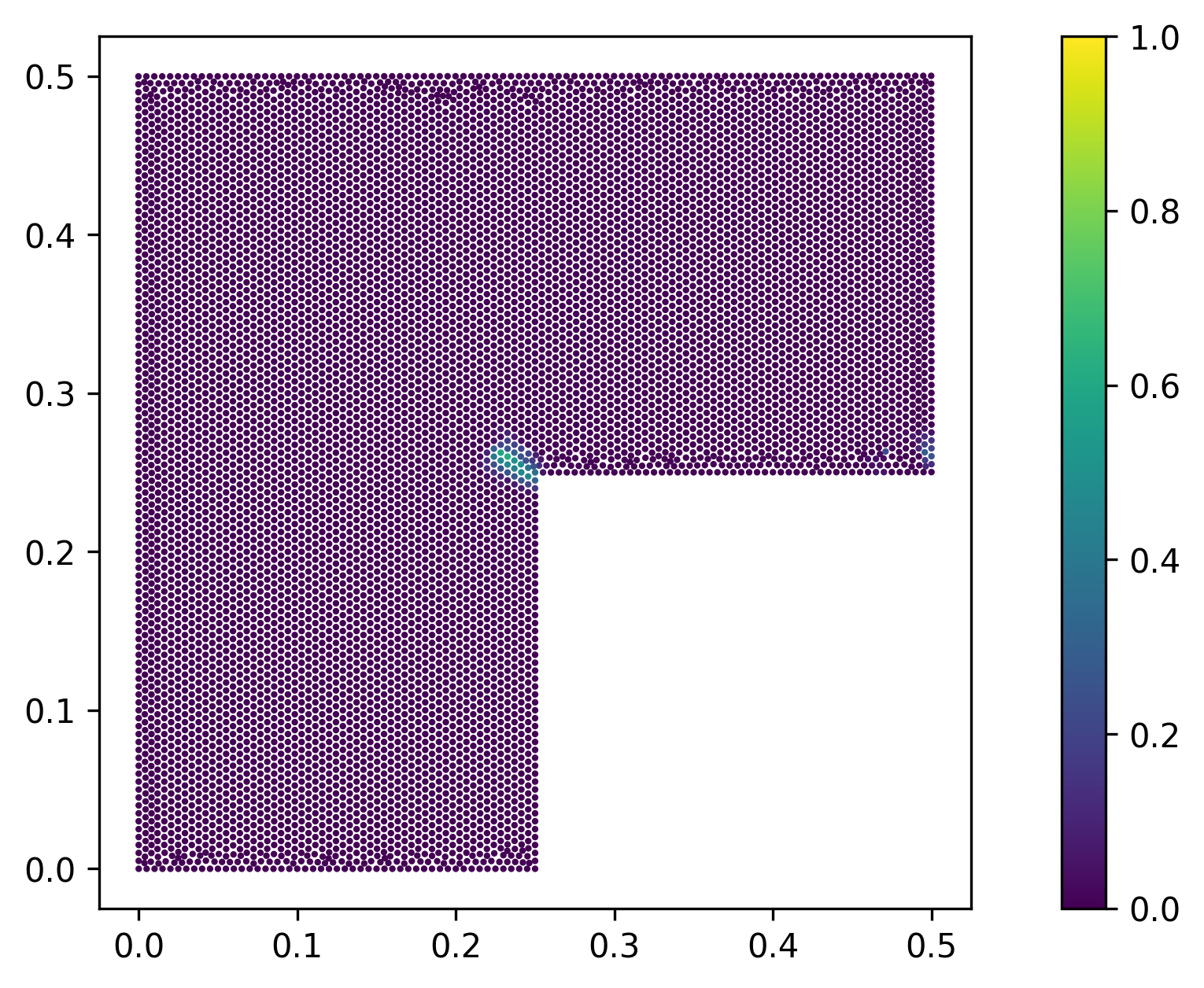}
\includegraphics[width=0.3\linewidth]{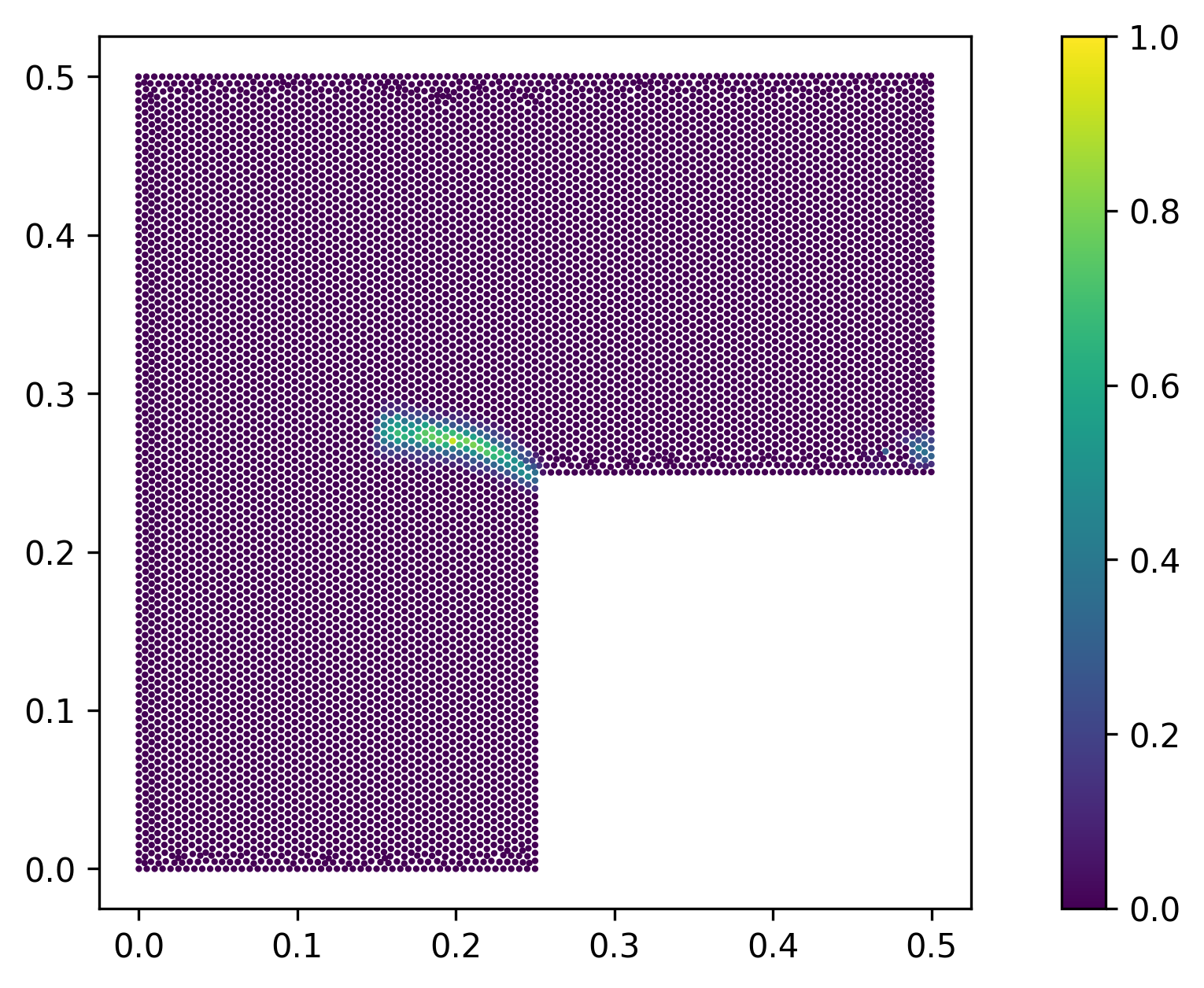}
\includegraphics[width=0.3\linewidth]{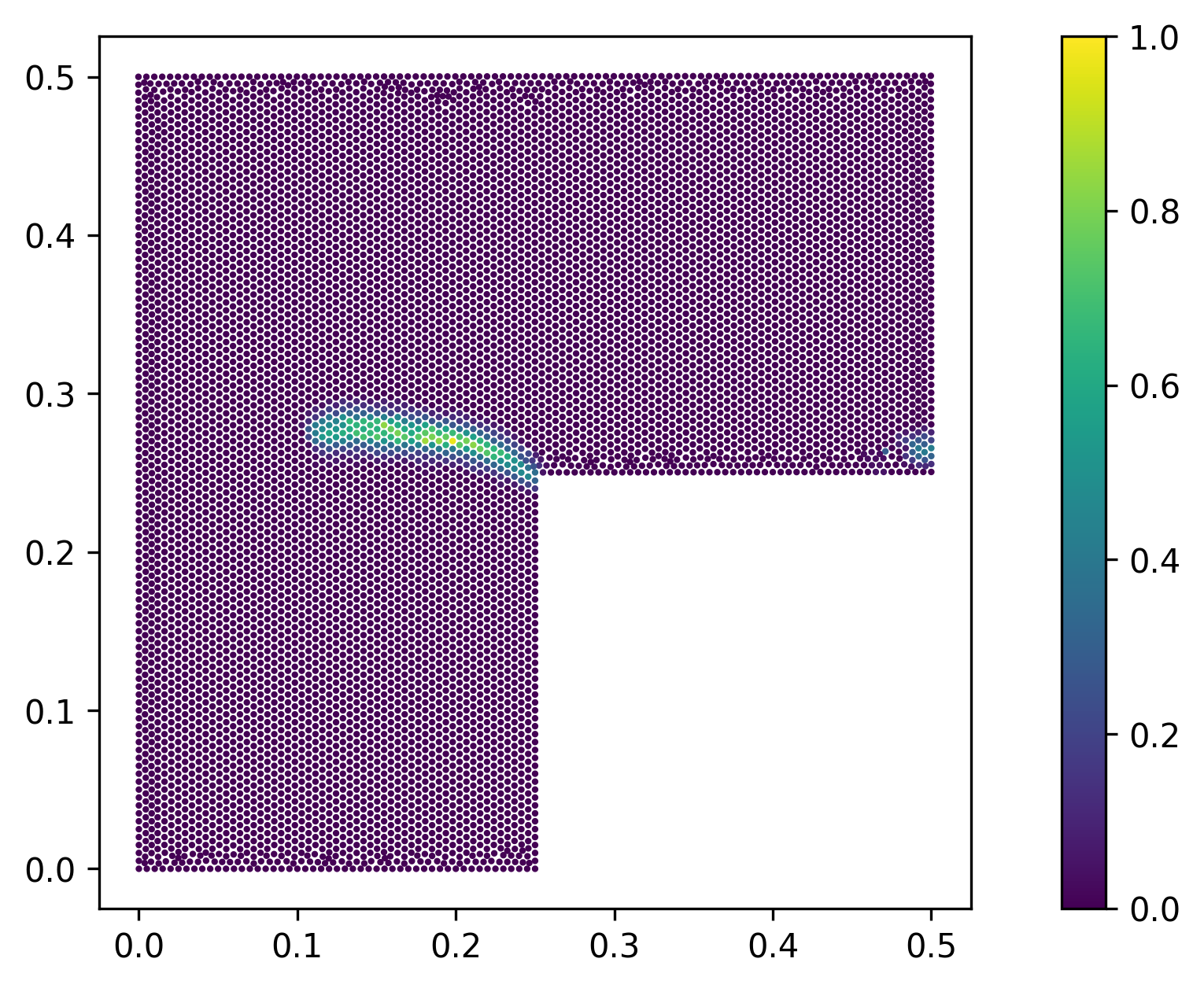}
\label{fig:7b}
}
		  \caption{Fracture in L-shaped domain}
		  \label{fig:L-shaped}
\end{figure}

We consider an L-shaped domain with geometry and boundary given in Fig. \ref{fig:diagram-L}.
The material parameters are chosen to be $\mu =$ 10.95 GPa and $G_c = $ 1000 J/m$^2$, which agree with the experiment given in \cite{Winkler}.
An unstructured grid with 8931 nodes is used to discretize the domain.
A constant peridynamic horizon of $\epsilon =$ 15 mm is used.
The bottom edge of the domain is clamped while a vertical displacement is applied in the bottom-right corner of the domain.
The initial displacement is taken to be $U^0 =$ 0.005 mm, and it is incremented by $\Delta U = $ 0.00076 mm at each load step.
Fig. \ref{fig:L-shaped} shows the bonds and the damage for the load steps $N= 30, 50$, and $60$.
The crack pattern is consistent with the crack pattern seen in experiment \cite{Winkler}.
The energy $E_N$ with respect to the load step $N$ is shown in Fig. \ref{fig:energy}.

\subsection{Double-notched tension test}
\label{sec:double-notched}
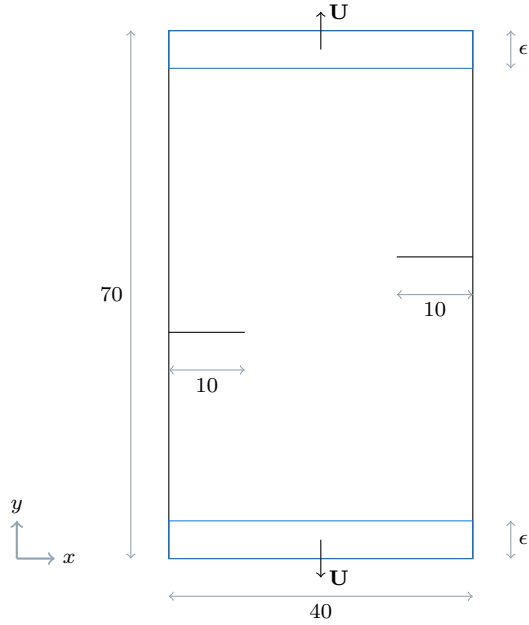
\begin{figure}[tb]
    \centering
    \begin{tikzpicture}
    \draw (0,0) rectangle (4,7);
    \draw  (0,3) -- (1,3);
    \draw  [<->,cadetgrey](0,2.5) -- (1,2.5);
    \node[below] at (0.5,2.5) {$10$};
    \draw  (3,4) -- (4,4);
    \draw  [<->, cadetgrey](3,3.5) -- (4,3.5);
    \node[below] at (3.5,3.5) {$10$};

    \draw[azure] (0,0) rectangle (4,0.5);
    \draw[<->,cadetgrey]  (4.5,0) -- (4.5,0.5);
    \node[right] at (4.5,0.25) {$\epsilon$};
    
    \draw[azure] (0,6.5) rectangle (4,7);
    \draw[<->,cadetgrey]  (4.5,6.5) -- (4.5,7);
    \node[right] at (4.5,6.75) {$\epsilon$};
    
    \draw[<->, cadetgrey] (0, -0.5) -- (4, -0.5);
    \node[below] at (2,-0.5) {$40$};
    
    \draw[<->, cadetgrey] (-0.5, 0) -- (-0.5, 7);
    \node[left] at (-0.5, 3.5) {$70$};
    
    \draw[->]  (2,0.25) -- (2,-0.25);
    \node[right] at (2,-0.25) {$\UU$};
    
    \draw[->]  (2,6.75) -- (2,7.25);
    \node[right] at (2,7.25) {$\UU$};
    
%
%

\draw[->,thick,cadetgrey] (-2,0) -- (-1.5,0);
\draw[<-,thick,cadetgrey] (-2,0.5) -- (-2,0);
\node[above] at (-2,0.5) {\small $y$};
\node[right] at (-1.5,0) {\small $x$};

    \end{tikzpicture}
    \caption{Double-notched domain under uniaxial tension loading}
    \label{fig:diagram-double-notched}
\end{figure}

\begin{figure}
\centering
\subfloat[Case 1: no offest: 0 mm]{
\includegraphics[width=0.3\linewidth]{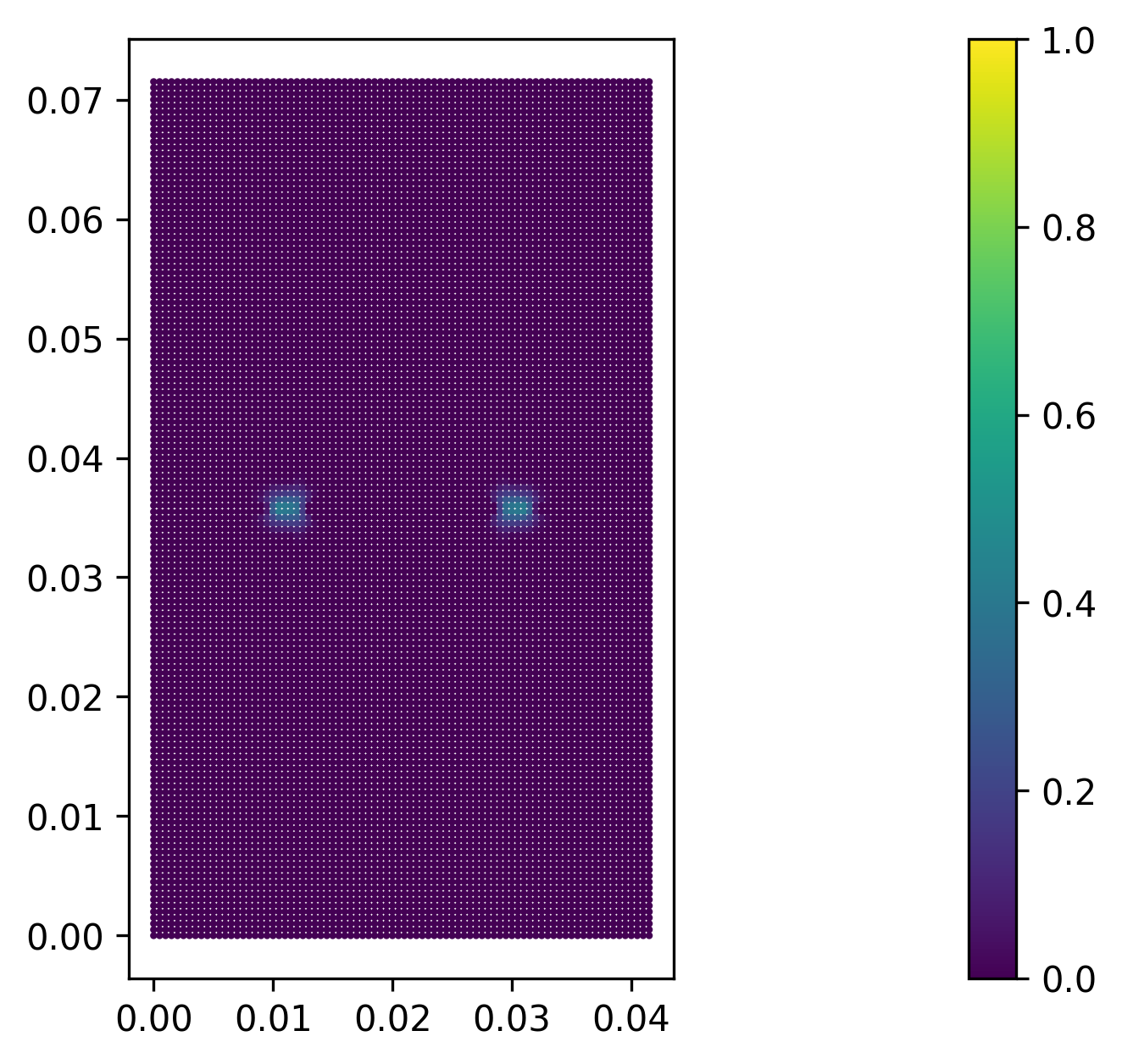}
\includegraphics[width=0.3\linewidth]{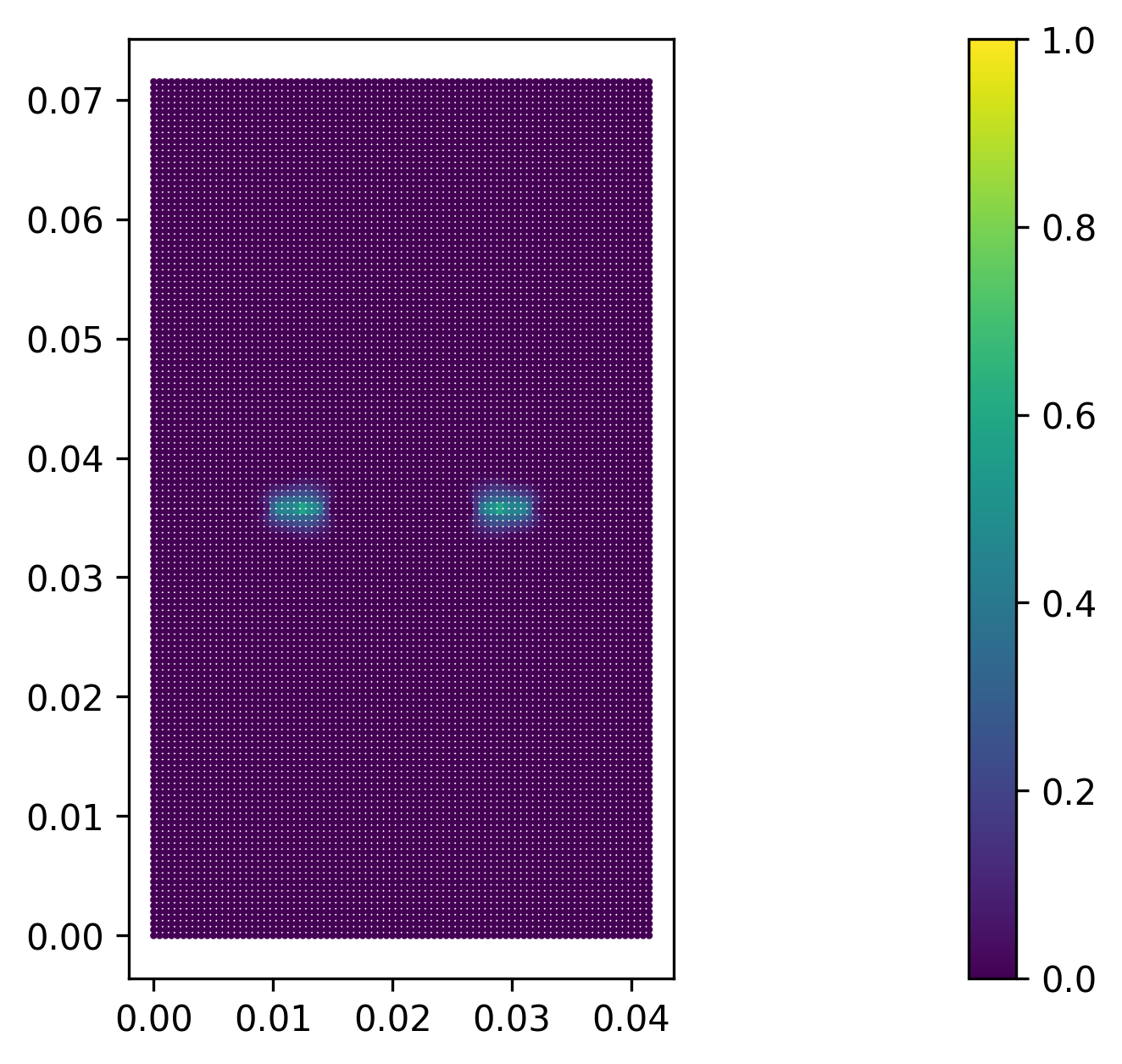}
\includegraphics[width=0.3\linewidth]{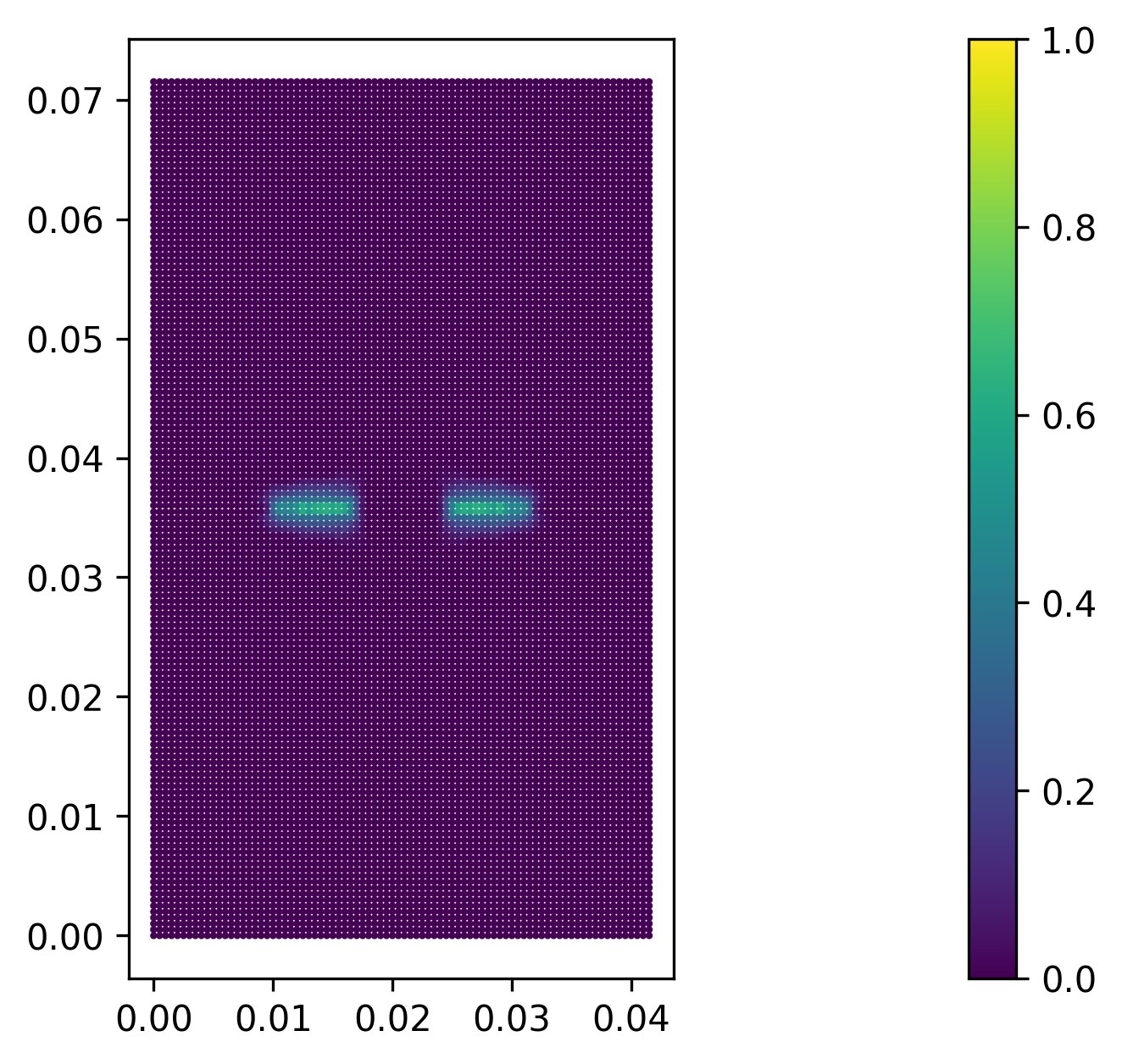}
}
\\
\subfloat[Case 2: small offest: 10 mm]{
\includegraphics[width=0.3\linewidth]{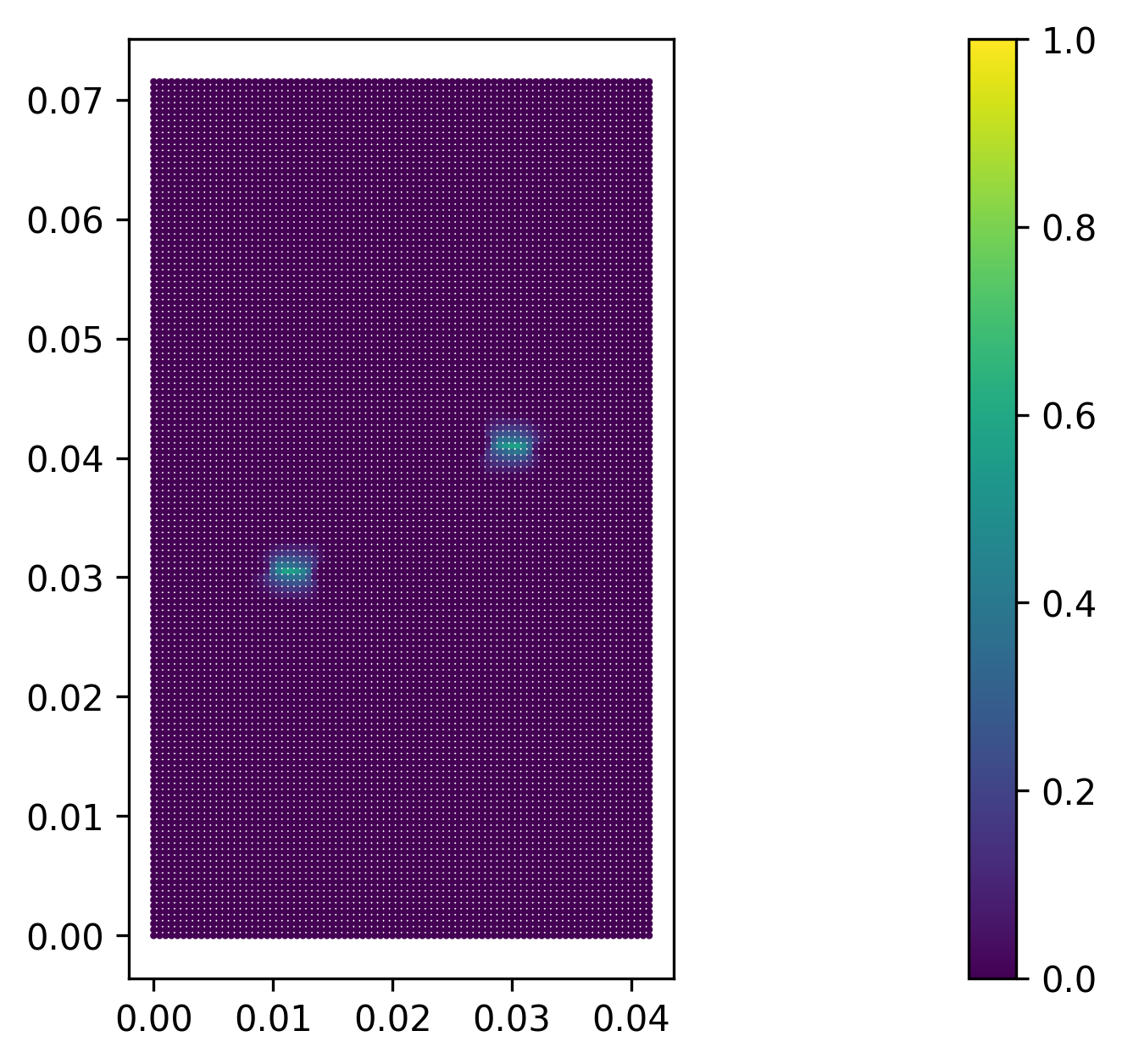}
\includegraphics[width=0.3\linewidth]{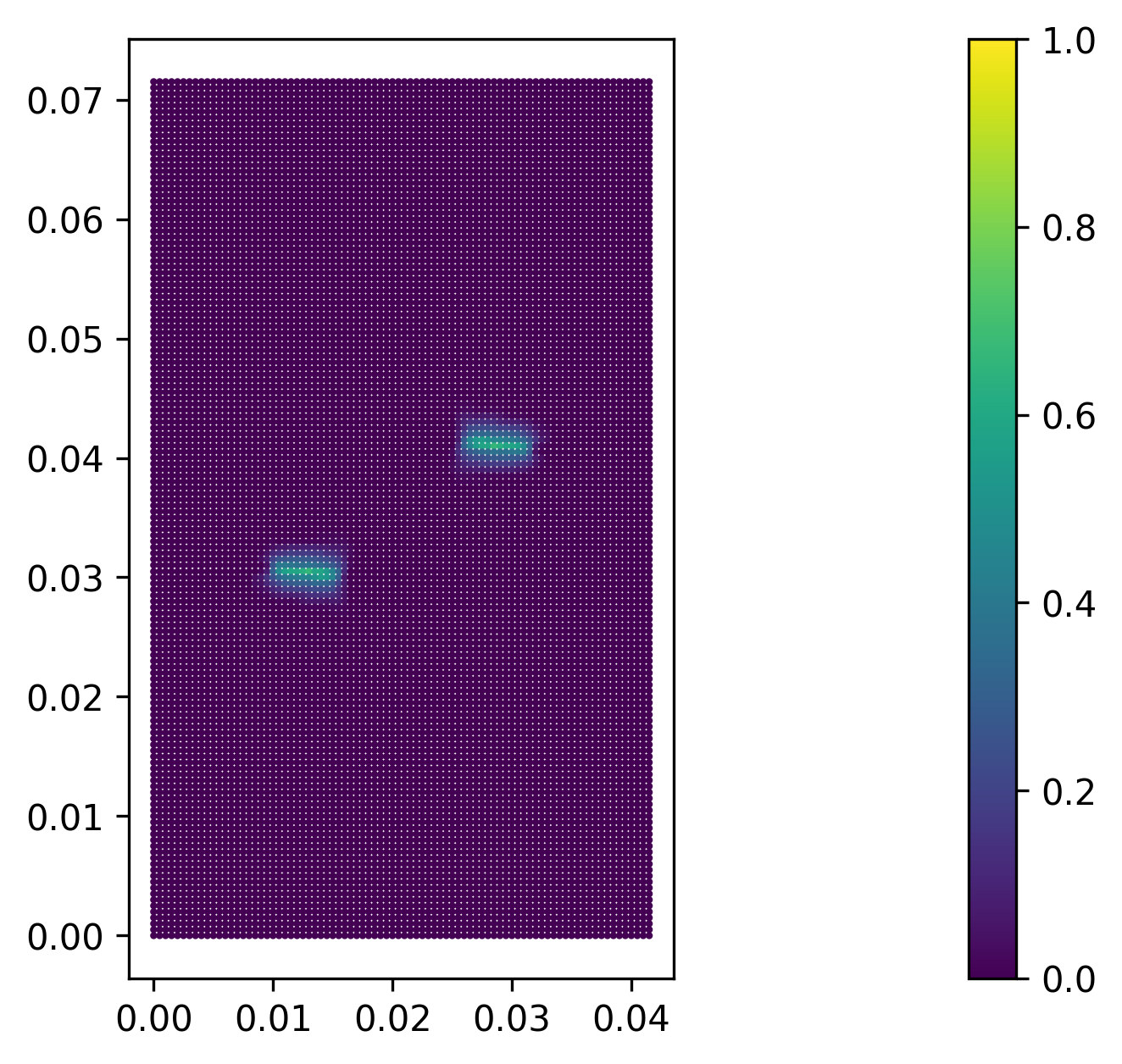}
\includegraphics[width=0.3\linewidth]{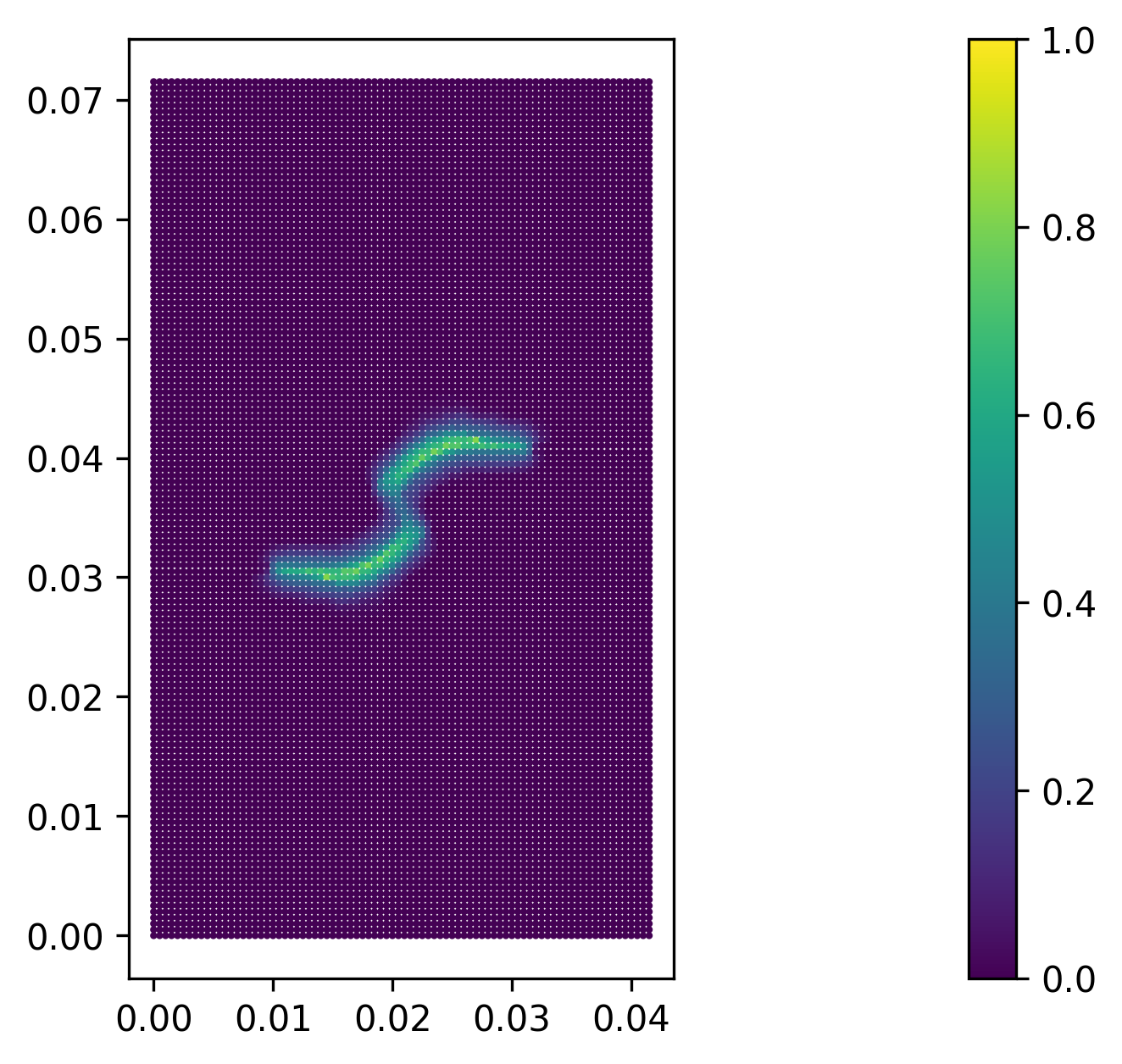}
}
\\
\subfloat[Case 3: large offset: 20 mm]{
\includegraphics[width=0.3\linewidth]{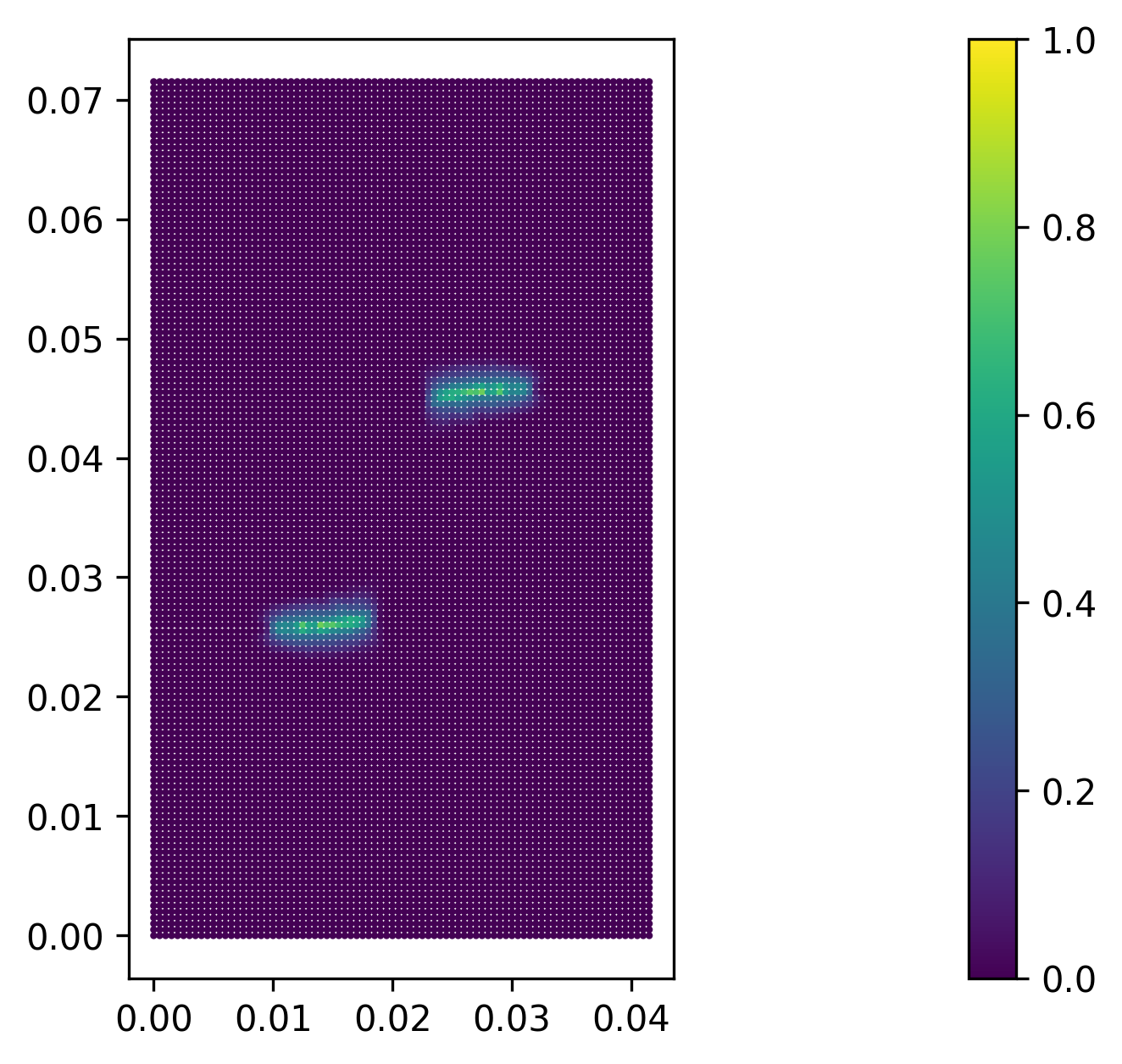}
\includegraphics[width=0.3\linewidth]{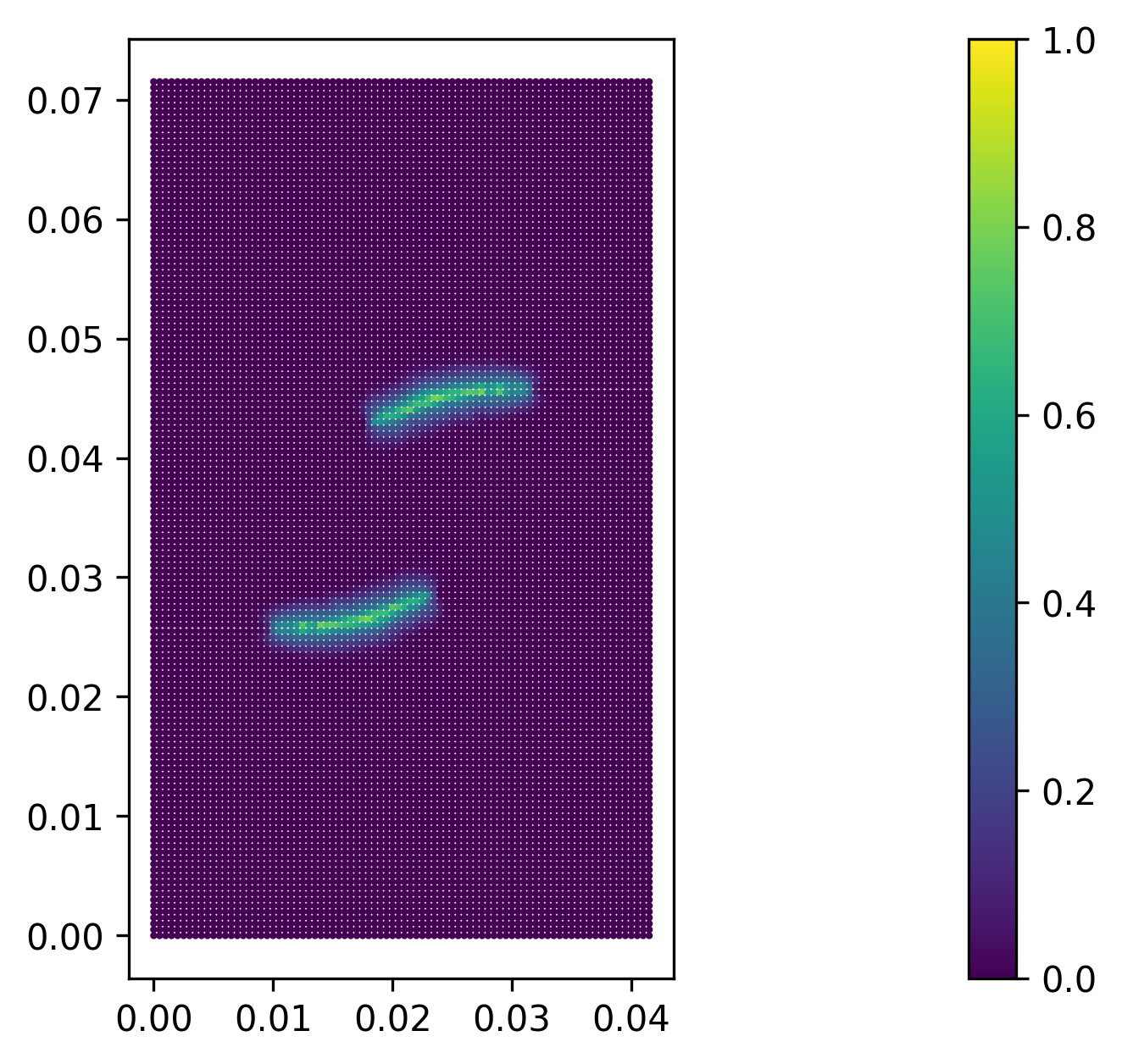}
\includegraphics[width=0.3\linewidth]{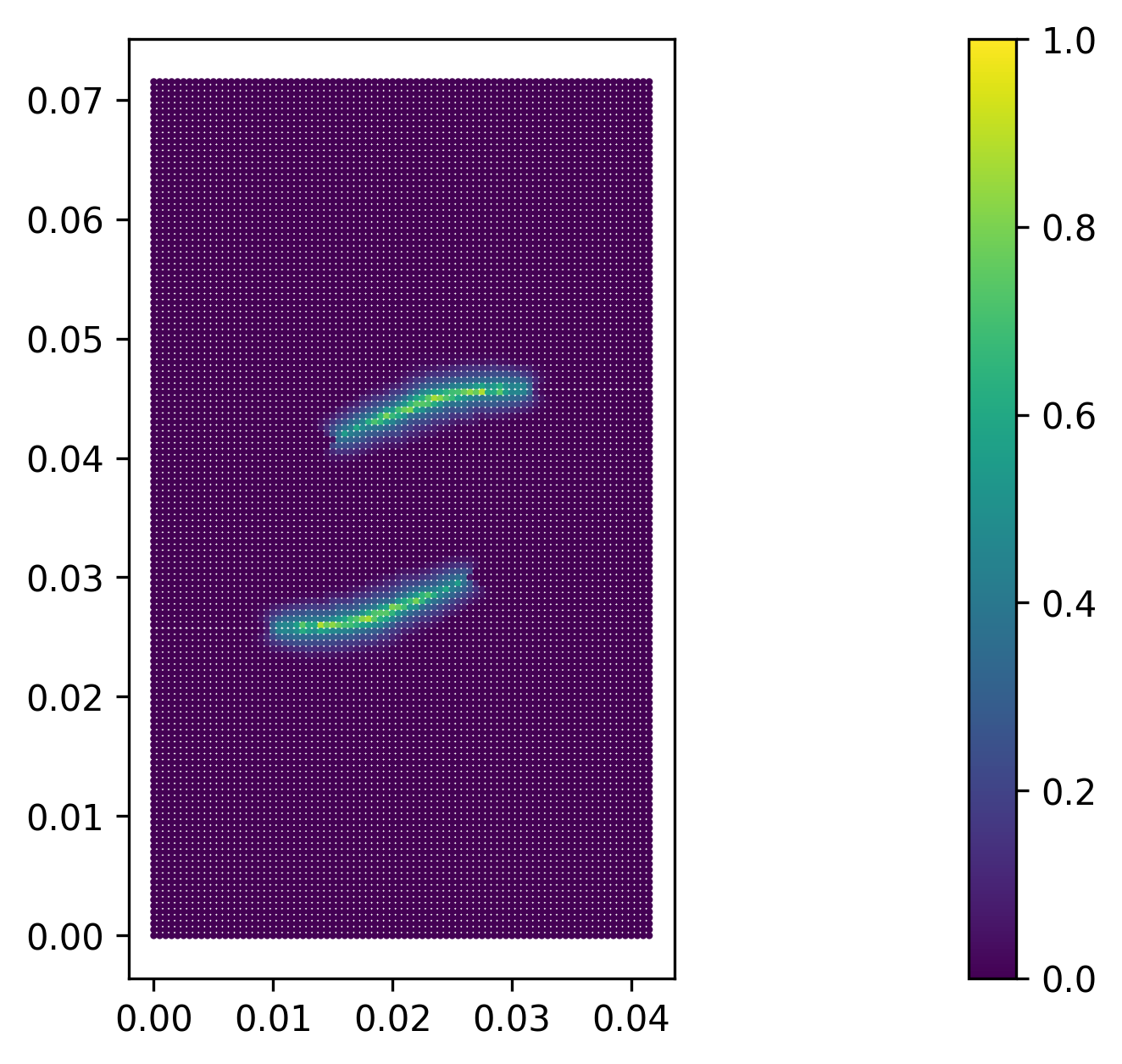}
}
		  \caption{Fracture patterns in double-notched tension test}
		  \label{fig:double-notched}
\end{figure}

We consider a double-notched rectangular domain given in Fig. \ref{fig:diagram-double-notched}, which was considered in \cite{zhao2018peridynamics}.
The material parameters are taken to be $E = $ 203 GPa and $G_c = 2700$ J/m$^2$.
The domain is a rectangle of dimension 40 mm $\times$ 70 mm. Two horizontal pre-notches of length 10 mm are present on the left and the right edges of the rectangle.
The vertical distance of the pre-notches are taken to be 0 mm, 10 mm, and 20 mm apart, respectively.
A regular grid is used and the mesh size is taken to be $h = 0.5$ mm. The peridynamic horizon is taken to be $\epsilon = m h$, where $m = 5$.
The domain is loaded under uniaxial tension in the vertical direction.
The initial applied displacement loading is taken to be $U^0 = 9.23 \times 10^{-5}$ mm. At each load step $N$, the displacement is increased by $\Delta U = 2.4 \times 10^{-5}$ mm.
The fracture patterns are shown in Fig. \ref{fig:double-notched}.
In the domain with pre-notches with zero vertical distance, the cracks grow inward in straight lines.
When the pre-notch distance is 10 mm, the cracks initially grow inward and eventually merge.
When the pre-notches are too far apart (20 mm), the cracks grow inward, bend slightly toward the center, but they do not merge.

\section{Conclusion}

In this paper a nonlocal quasistatic model is developed for the evolution of interacting cracks. Using this model, numerical examples are presented. The approach is implicit and is based on fixed point methods. For each load step it is seen that  the  Newton convergence of the residual as measured by the maximum norm is fast. For the straight crack it takes at most four iterations in the presence of crack growth before the residual lies below a tolerance of $10^{-5}$.
We have proved that the fracture evolution decreases stored elastic energy of  the intact material with each load step as the cracks advance. This holds true theoretically provided the load increments are chosen sufficiently small. This is borne out in the numerical examples. 
The numerical examples show that crack patterns emerge from the field theory in the quasistatic context. 
 
%

\section*{Acknowledgments}
This material is based upon work supported by the U. S. Army Research Laboratory and the U. S. Army Research Office under Contract/Grant Number W911NF-19-1-0245.

\section*{Appendix}

We provide the proof of Proposition \ref{initialize}  below.
Substitution of  $\uu^0=0$ into \eqref{eq: hessianN} gives
\begin{align}
      \label{eq: hessian0}
   \LL'_0[\uu^0]\Delta\uu(\xx) = 
   -\int\limits_{\hat{H}_\epsilon(\xx) \cap \Omega} \frac{J^\epsilon(\abs{\yy - \xx} )}{\epsilon^{n+1} \omega_n} g''\left(0\right)S(\yy, \xx, \Delta \uu)  \ee_{\yy - \xx} d\yy,
\end{align}	
and we directly  verify as in \cite{BhattacharyaLipton} that $Ker\{\LL'[\uu^0]\}=0$ and $\mathbb{A}[\uu^0]-\gamma\mathbb{I}>0$ for some $\gamma>0$. So $\LL'[\uu^0]^{-1}$ exists and is a bounded linear operator mapping $\mathcal{V}$ onto $\mathcal{V}$.

Now proceeding as in Section \ref{energy inequality} we conclude there is closed ball surrounding $\uu^0$ given by $\overline{B(R,\uu^0)}=\{\uu:\, ||\uu-\uu^0||_\infty\leq R\}$ of radius $R>0$ and center  for which $\LL'_N[\uu]^{-1}$ is well defined. Consequently  there is a positive constant $K_0>0$  for which $K_0||\ww||_\infty\leq ||\LL'[\uu]\ww||_\infty$ for any fixed $\uu$ in $\overline{B(R,\uu^0)}$ and for all $\ww\in\mathcal{V}$.  With this in hand, we  now show that $\uu^0=0$ is the only solution of the boundary value problem  among all functions  in  $\overline{B(R,\uu^0)}$. To see this suppose there is another solution $\uu$ in  $\overline{B(R,\uu^0)}$ to the boundary value problem $\LL[\uu^0](\xx)=0$,  for $\xx$ in $D$ and $\uu=0$ on $\Omega_d$. Applying the fundamental theorem of calculus and the mean value theorem with  gives a $0\leq \tilde{t}\leq 1$ for which
\begin{align}
\label{unique}
0=||\LL_0[\uu]-\LL_0[\uu^0]||=||\int_0^1\LL'[t\uu]\uu\,dt||=||\LL'[\tilde{t}\uu]\uu||\geq K_0||\uu||,
\end{align}
so $\uu=0$. Hence $\uu^0=0$ is the only solution in this ball.


%
%

\bibliographystyle{spmpsci}      
\bibliography{references}   

\end{document}